\documentclass[12pt,thmsa]{article}
\usepackage{amsfonts}
\usepackage[dvips]{graphics}

\def\mylabel#1{\label{#1}}
\newtheorem{theorem}{Theorem}[section]
\newtheorem{lemma}[theorem]{Lemma}
\newtheorem{corollary}[theorem]{Corollary}
\newtheorem{proposition}[theorem]{Proposition}
\newtheorem{example}[theorem]{Example}
\newtheorem{remark}[theorem]{Remark}

\newtheorem{hypothesis}[theorem]{Hypothesis}

\def\bit{\begin{itemize}}
\def\eit{\end{itemize}}
\reversemarginpar   
\def\bc{\begin{center}}
\def\ec{\end{center}}
\def\bthm{\begin{theorem}}
\def\ethm{\end{theorem}}
\def\bcor{\begin{corollary}}
\def\ecor{\end{corollary}}
\def\bprop{\begin{proposition}}
\def\eprop{\end{proposition}}
\def\blem{\begin{lemma}}
\def\elem{\end{lemma}}
\def\bex{\begin{example} {\rm }
\def\eex{\end{example} }}
\def\brem{\begin{remark}}
\def\erem{\end{remark}}
\def\prf{\noindent{\bf Proof~: }}
\def\bdes{\begin{description}}
\def\edes{\end{description}}

\def\iti{\item[(i)]}
\def\itii{\item[(ii)]}
\def\itiii{\item[(iii)]}
\def\itiv{\item[(iv)]}

\def\beq{\begin{equation}}
\def\eeq{\end{equation}}

\def\ben{\begin{enumerate}}
\def\een{\end{enumerate}}

\def\beqar{\begin{eqnarray}}
\def\eeqar{\end{eqnarray}}
\def\beqarr{\begin{eqnarray*}}
\def\eeqarr{\end{eqnarray*}}
\def \non{{\nonumber}}

\def\RR{{\mathbb R}}  

\def\cA{\mathcal{A}} \def\cB{\mathcal{B}} \def\cC{\mathcal{C}}
  \def\cF{\mathcal{F}}
  
 \def\cK{\mathcal{K}} \def\cL{\mathcal{L}}
\def\cM{\mathcal{M}}  
\def\cP{\mathcal{P}}  
  
  \def\cX{\mathcal{X}}

\def\Rp{{\mathbb R}_+}   

\def\qed{\hspace{.1in}{\bf QED}}

\def\P{{\mathsf P}} 
\def\Q{{\mathsf Q}} 

\def\E{{\mathsf E}} 

\def\NN{{\mathbb N}}       

\def\one{{\bf 1}}

\def\eps{\epsilon}

\def\p{\varphi}

\def\la{\langle}

\def\ra{\rangle}

\def\part{\partial}

\def\d#1dt{\frac{d#1}{dt}}    

\begin{document}
\title{Self-interacting diffusions IV: Rate of convergence\thanks{We acknowledge financial support from the  Swiss National Science Foundation Grant 200021-103625/1}}
\author{{\bf Michel Bena\"{\i}m}\\
Universit\'e de Neuch\^atel, Suisse \and
{\bf Olivier Raimond}\\Universit\'e Paris Ouest Nanterre la d\'efense, France}

\maketitle

\begin{center} {\large\bf Abstract} \end{center}
Self-interacting diffusions are processes living on a compact Riemannian manifold defined by a stochastic differential equation with a drift term depending on the past empirical measure $\mu_t$ of the process. The asymptotics of $\mu_t$ is governed by a deterministic dynamical system
and under certain conditions $(\mu_t)$ converges almost surely towards a deterministic measure $\mu^*$ (see Bena\"{\i}m, Ledoux, Raimond (2002) and Bena\"{\i}m, Raimond (2005)). We are interested here in the rate of convergence of $\mu_t$ towards $\mu^*$. A central limit theorem is proved. In particular, this shows that greater is the interaction repelling faster is the convergence.

\section{Introduction}

\medskip
\subsection*{Self-interacting diffusions}

Let $M$ be a smooth compact Riemannian manifold and  $V:M\times M\to \RR$  a sufficiently smooth mapping\footnote{The mapping $V_x:M\to\RR$ defined by $V_x(y)=V(x,y)$ is $C^2$ and its derivatives are continuous in $(x,y)$}.
For all finite Borel measure $\mu$, let  $V\mu : M \to \RR$ be the smooth function defined by
$$V\mu(x)=\int_M V(x,y)\mu(dy).$$
Let $(e_\alpha)$ be a finite  family of vector fields on $M$  such that $$\sum_\alpha e_\alpha(e_\alpha f)(x)=\Delta f(x),$$ where $\Delta$ is the Laplace operator on $M$ and $e_{\alpha}(f)$ stands for the Lie derivative of $f$ along $e_{\alpha}.$
Let $(B^\alpha)$ be a family of independent Brownian motions.

A {\em self-interacting diffusion}  on $M$ associated to $V$ can be defined as the solution to the stochastic differential equation (SDE)
$$dX_t=\sum_\alpha e_\alpha(X_t)\circ dB^\alpha_t - \nabla (V\mu_t) (X_t) dt.$$
where $$\mu_t=\frac{1}{t}\int_0^t\delta_{X_s}ds$$ is the empirical occupation measure of $(X_t).$

In absence of drift (i.e~$V = 0$), $(X_t)$ is just a Brownian motion on $M$
but in general  it defines a non Markovian process whose behavior at time $t$  depends on its past trajectories through $\mu_t.$
This type of process was introduced in  Benaim, Ledoux and Raimond (2002) (hence after referred as \cite{blr}) and further analyzed in a series of papers by Benaim and Raimond (2003, 2005, 2007) (hence after referred as \cite{br1}, \cite{br2} and \cite{br3}).
We refer the reader to these papers for more details and especially to \cite{blr} for a detailed construction of the process and its elementary properties.
For a general overview of processes with reinforcement we refer the reader to the recent survey paper by Pemantle (2007) (\cite{Pem07}).

\subsection*{Notation and Background}

 \paragraph{Standing Notation}
We let $\cM(M)$ denote the space of finite Borel measures on $M,$ $\cP(M) \subset \cM(M)$ the space of probability measures. If $I$ is a metric space (typically, $I = M, \RR^+ \times M$ or $[0,T] \times M$) we let   $C(I)$ denote the space of real valued continuous functions on $I$ equipped with the topology of uniform convergence on compact sets. When $I$ is compact and  $f \in C(I)$ we let
$\|f\| = \sup_{x \in I} |f(x)|.$
 The normalized Riemann measure on $M$ will be denoted by $\lambda.$

Let $\mu \in \cP(M)$ and $f : M \to \RR$ a nonnegative or $\mu-$integrable Borel function. We  write $\mu f$ for $\int f d\mu,$ and $f \mu$ for the measure defined as $f \mu(A) = \int_A f d\mu.$   We let $L^2(\mu)$ denote the space of such functions  for which $\mu |f|^2 < \infty,$ equipped with the inner product $$\langle f, g \rangle_{\mu} = \mu (fg)$$
and the norm
$$ \|f\|_\mu = \sqrt{\mu f^2}.$$ We simply write $L^2$ for $L^2(\lambda).$

Of fundamental importance in the analysis of the asymptotics of $(\mu_t)$ is the mapping
 $\Pi  : {\cal M}(M) \to {\cal P}(M)$  defined by
\beq \label{eq:PI}
\Pi(\mu)  = \xi( V{\mu}) \lambda
\eeq
where $\xi : C(M) \to C(M)$ is the function defined by
\beq \label{eq:defxi}
\xi(f)(x) = \frac{e^{-f(x)}}{\int_M e^{-f(y)}\lambda(dy)}.
\eeq
In \cite{blr}, it is shown that the asymptotics of $\mu_t$ can be precisely  related to the long term behavior of a certain semiflow on  $\cP(M)$ induced by  the ordinary differential equation (ODE) on $\cM(M):$
\beq \dot{\mu}=-\mu+\Pi(\mu). \label{basicode} \eeq
Depending on the nature of $V,$ the dynamics of (\ref{basicode}) can either be convergent or nonconvergent  leading to  similar behaviors for $\{\mu_t\}$(see \cite{blr}).  When $V$ is symmetric, (\ref{basicode}) happens to be   a {\em quasigradient} and the following convergence result hold.
 \bthm [\cite{br2}]
\label{th:brmain}
Assume that $V$ is symmetric. i.e. $V(x,y) = V(y,x).$ Then
the limit set of $\{\mu_t\}$ (for the
 topology of weak* convergence)  is almost surely a compact connected subset of
 $$\mathsf{Fix}(\Pi) = \{ \mu \in {\cal P}(M) : \mu = \Pi(\mu)\}.$$
\ethm
\smallskip
In particular, if $\mathsf{Fix}(\Pi)$ is finite then $(\mu_t)$
converges almost surely toward a fixed point of $\Pi$.  This holds for a generic function $V$ (see \cite{br2}).

Sufficient conditions ensuring that $\mathsf{Fix}(\Pi)$ has cardinal one  are as follows:
\bthm [\cite{br2}, \cite{br3}]
\mylabel{th:globconv}
Assume that $V$ is symmetric and that one of the two following conditions hold
\bdes
 \iti Up to an additive constant $V$ is a Mercer kernel,
That is $$V(x,y) = K(x,y) + C$$
and $$\int K(x,y) f(x) f(y) \lambda(dx) \lambda(dy) \geq 0 $$
for all $f \in L^2.$
\itii  For all $x \in M, y \in M, u \in T_x M, v \in T_y M$
$$\mathsf{Ric}_x(u,u) + \mathsf{Ric}_y(v,v) +  \mathsf{Hess}_{x,y}
V((u,v),(u,v)) \geq K (\|u\|^2 + \|v\|^2)$$
where $K$ is some positive constant. Here $\mathsf{Ric}_x$ stands for the Ricci tensor at $x$ and $ \mathsf{Hess}_{x,y}$ is the Hessian of $V$ at $(x,y).$
\edes
Then  $\mathsf{Fix}(\Pi)$ reduces to a singleton $\{\mu^*\}$ and
$\mu_t \to \mu^*$ with probability one.
\ethm
As observed in \cite{br3} the condition $(i)$ in Theorem \ref{th:globconv} seems well suited to describe {\em self-repelling diffusions}. On the other hand, it is not  clearly related to the geometry of $M.$ Condition $(ii)$
has a more geometrical  flavor and is robust to smooth perturbations (of $M$ and $V$). It can be seen as a Bakry-Emery type condition for self interacting diffusions.

In \cite{br2}, it is also proved that every stable (for the ODE (\ref{basicode})) fixed point of $\Pi$ has a positive probability to be a limit point for $\mu_t$; and any unstable fixed point cannot be a limit point for $\mu_t$.

\subsection*{Organisation of the paper}
Let $\mu^*\in \mathsf{Fix}(\Pi)$. We will assume that
\begin{hypothesis}\label{hyp1}
$\mu_t$ converges a.s. towards $\mu^*$.
\end{hypothesis}
Sufficient conditions are given by Theorem \ref{th:globconv}

In this paper we intend to study the rate of this convergence. Let
$$\Delta_t = e^{t/2}(\mu_{e^t} - \mu^*).$$ It will be shown  that, under some conditions to be specified later,  for all  $g = (g_1, \ldots, g_n) \in C(M)^n$ the process
 $$\left [  \Delta_s g_1, \ldots, \Delta_s g_n,V\Delta_{s} \right ]_{s \geq t} $$ converges in law, as $t \to \infty,$ toward a certain  stationary Ornstein-Uhlenbeck process $(Z^g,Z)$ on
$\RR^n \times C(M).$ This process is defined in Section \ref{sec:OU1}. The main result is stated in section \ref{tlcmu} and  some examples are  developed.  It is in particular observed that a strong repelling interaction gives a faster convergence.
The section \ref{Vmu} is a proof section. The appendix, section \ref{sec:app}, contains general material on random  variables and Ornstein-Uhlenbeck processes on $C(M).$

\medskip In the following $K$ (respectively $C$)
denotes a positive constant (respectively a positive random
constant). These constants may change from line to line.

\section{The Ornstein-Uhlenbeck process $(Z^g,Z).$}
\label{sec:OU1}

\label{gmu}
Throughout all this section we let $\mu\in \cP(M)$.
For $x \in M$ we set $V_x:M\to\RR$ defined by $V_x(y) = V(x,y).$
\subsection{The operator $G_\mu$}
Let $g \in C(M)$ and
let $G_{\mu,g}:\RR\times C(M)\to \RR$ be the linear operator defined by
\beq G_{\mu,g} (u,f)=u/2+\hbox{Cov}_\mu (g,f),\eeq
where $\hbox{Cov}_\mu$ is the covariance on $L^2(\mu)$, that is  the bilinear form acting on $L^2\times L^2$ defined by
$$\hbox{Cov}_\mu(f,g)=\mu (fg) - (\mu f)(\mu g).$$
We define the linear operator $G_{\mu}: C(M)\to C(M)$ by
\beqar
\label{eq:Gmu} G_\mu f(x) &=& G_{\mu,V_x}(f(x),f)\\
&=& f(x)/2 + \hbox{Cov}_\mu(V_x,f).\non\eeqar
It is easily seen that $\|G_{\mu} f\| \leq (2 \|V\| + 1/2)\|f\|.$ In particular, $G_{\mu}$ is a  bounded operator.
Let $\{e^{-t G_{\mu}}\}$ denotes the semigroup acting on $C(M)$ with generator  $-G_{\mu}.$
From now on we will assume the following:
\begin{hypothesis} \label{hyp:mu}
There exists $\kappa > 0$ and $\widehat{\lambda}\in\cP(M)$ such that $\mu<<\widehat{\lambda}$ with $\|\frac{d\mu}{d\widehat{\lambda}}\|_\infty<\infty$, $\lambda$ and $\widehat{\lambda}$ are equivalent measures with $\|\frac{d\lambda}{d\widehat{\lambda}}\|_\infty<\infty$ and $\|\frac{d\widehat{\lambda}}{d\lambda}\|_\infty<\infty$,  and such that for all $f\in L^2(\widehat{\lambda})$,
$$\langle G_\mu f,f\rangle_{\widehat{\lambda}}\geq \kappa \|f\|^2_{\widehat{\lambda}}.$$
\end{hypothesis}

Let
$$\lambda(-G_\mu)=\lim_{t \to \infty} \frac{\log(\|e^{-t G_{\mu}}\|)}{t}.$$
This limit exists by subadditivity. Then
\blem \label{lemhyp} Hypothesis \ref{hyp:mu} implies that
$\lambda(-G_\mu)\le -\kappa<0$. \elem

\prf For all $f\in L^2(\widehat{\lambda})$, \beqarr
\frac{d}{dt}\|e^{-tG_\mu}f\|_{\widehat{\lambda}}^2 &=& -2\langle G_\mu
e^{-tG_\mu}f, e^{-tG_\mu}f \rangle_{\widehat{\lambda}}\\
&\leq& -2\kappa \|e^{-tG_\mu}f\|_{\widehat{\lambda}}. \eeqarr This implies that
$\|e^{-tG_\mu}f\|_{\widehat{\lambda}}\leq e^{-\kappa t}\|f\|_{\widehat{\lambda}}$.

Denote by $g_t$ the solution of the differential equation
$$\frac{dg_t}{dt}=\hbox{Cov}_\mu(V_x,g_t)$$
with $g_0=f$, where $f\in C(M)$. Note that $e^{-tG_\mu}f=e^{-t/2}g_t$. It is straightforward to check that (using the fact that $\|\frac{d\mu}{d{\widehat{\lambda}}}\|_\infty <\infty$)
$$\frac{d}{dt}\|g_t\|_{\widehat{\lambda}}\leq K \|g_t\|_{\widehat{\lambda}}$$
with $K$ a constant depending only on $V$ and $\mu$. Thus $$\sup_{t\in [0,1]}\|g_t\|_{\widehat{\lambda}} \leq K \|f\|_{\widehat{\lambda}}.$$
Now, since for all $x\in M$ and $t\in [0,1]$
$$\left|\frac{d}{dt}g_t(x)\right|\leq K \|g_t\|_{\widehat{\lambda}} \leq K \|f\|_{\widehat{\lambda}},$$
we have $\|g_1\|\leq K \|f\|_{\widehat{\lambda}}$. This implies that
$$\|e^{-G_\mu}f\|\leq K \|f\|_{\widehat{\lambda}}.$$

Now for all $t>1$, and $f\in C(M)$,
\beqarr
\|e^{-tG_\mu}f\|
&=& \|e^{-G_\mu}e^{-(t-1)G_\mu}f\|\\
&\leq& K\|e^{-(t-1)G_\mu}f\|_{\widehat{\lambda}}\\
&\leq& Ke^{-\kappa(t-1)}\|f\|_{\widehat{\lambda}}\\
&\leq& Ke^{-\kappa t}\|f\|_\infty.
\eeqarr
This implies that $\|e^{-tG_\mu}\|\leq Ke^{-\kappa t}$, which proves the lemma.
\qed

\medskip
The {\em adjoint} of $G_{\mu}$ is the operator on $\cM(M)$ defined by the relation
$$m (G_{\mu} f) = (G_{\mu}^* m) f$$ for all $m \in \cM(M)$ and $f \in C(M).$
It is not hard to verify that
\beq
\label{eq:Gmu*}
G_{\mu}^* m   = \frac{1}{2}m + (Vm) \mu - (\mu (Vm)) \mu.
\eeq

\subsection{The generator $A_\mu$ and its inverse $\Q_\mu$}
Let $H^2$ be the Sobolev space of real valued functions on $M$, associated with the norm $\|f\|^2_H=\|f\|_\lambda^2+\|\nabla f\|_\lambda^2$.
Since $\Pi(\mu)$ and $\lambda$ are equivalent measures with continuous Radon-Nykodim derivative, $L^2(\Pi(\mu))=L^2(\lambda):=L^2$.
We denote by $K_\mu$ the projection operator, acting on $L^2(\Pi(\mu))$,  defined by
$$K_\mu f=f-\Pi(\mu)f.$$
We denote by $A_\mu$ the operator acting on $H^2$ defined by
$$A_\mu f=\frac{1}{2}\Delta f -\langle \nabla V\mu,\nabla f\rangle.$$
Note that for $f$ and $g$ in $L^2$,
$$\langle A_\mu f,g\rangle_{\Pi(\mu)}=-\frac{1}{2}\int \langle \nabla f,\nabla g\rangle(x) \Pi(\mu)(dx)$$
where $\langle \cdot,\cdot \rangle$ denotes the Riemannian inner product on $M.$

For all $f\in C(M)$ there exists $\Q_\mu f\in H^2$ such that $\Pi(\mu)(\Q_\mu f)=0$ and
\beq f-\Pi(\mu)f= K_\mu f = - A_\mu \Q_\mu f.\eeq
Note that if $P^\mu_t$ denotes the semigroup with generator $A_\mu$, then
$$\Q_\mu f = \int_0^\infty P^\mu_t K_\mu f dt.$$
Since there exists $p^\mu_t(\cdot,\cdot)$ such that
$$P^\mu_tf(x) = \int_M p_t^\mu(x,y) f(y) \Pi(\mu)(dy),$$
we have
$$\Q_\mu f(x)=\int_M q_\mu(x,y)f(y)\Pi(\mu)(dy)$$
where
$$q_\mu(x,y)=\int_0^\infty (p_t^\mu(x,y)-1) dt.$$

Then, as  shown in \cite{blr}, $\Q_\mu f$ is $C^1$
and there exists a constant $K$ such
that for all $f\in C(M)$ and $\mu\in\cP(M)$,
\beqar \|\Q_\mu f\|_\infty\leq K\|f\|_\infty\label{sgest1}\\
\|\nabla\Q_\mu f\|_\infty\leq K\|f\|_\infty.\label{sgest2}
\eeqar

\medskip
Finally, note that for $f$ and $g$ in $L^2$,
\begin{eqnarray}
\int\langle \nabla \Q_\mu f,\nabla \Q_\mu g\rangle(x)\Pi(\mu)(dx)
&=& -2\langle A_\mu\Q_\mu f,\Q_\mu g\rangle_{\Pi(\mu)} \label{Chat}\\
&=& 2\langle f,\Q_\mu g\rangle_{\Pi(\mu)}. \non
\end{eqnarray}

\subsection{The covariance $C_\mu$}
We let $\widehat{C}_\mu$ denote the bilinear continuous form $\widehat{C}_\mu:C(M)\times C(M)\to\RR$ defined
by $$\widehat{C}_\mu(f,g)= 2\langle f,\Q_\mu g\rangle_{\Pi(\mu)} .$$
This form is symmetric (see its expression given by (\ref{Chat})). Note also that for some constant depending on $\mu$,
$$|\widehat{C}_\mu(f,g)|\leq K \|f\|\times\|g\|.$$

We let $C_\mu$ denote the mapping $C_\mu:M\times M\to\RR$ defined
by $$C_\mu(x,y)=\widehat{C}_\mu(V_x,V_y).$$ Then $C_\mu$ is a
covariance function (or a Mercer kernel), i.e. it is continuous,
symmetric and $\sum_{i,j}\lambda_i\lambda_j C_{\mu}(x_i,x_j)\geq 0$.

\subsection{The process $Z$}\label{secoucmu}
We now  define an Ornstein-Uhlenbeck process on $C(M)$ of covariance $C_{\mu}$ and drift $-G_{\mu}.$ This heavily relies on the general construction given in the appendix.

A  {\em Brownian motion on $C(M)$ with covariance $C_{\mu}$} is a $C(M)$-valued stochastic process  $W  = \{ W_t\}_{ t \geq 0}$ such that
\bdes
\iti  $W_0 = 0;$
\itii $t \mapsto W_t$ is continuous;
\itiii For every finite subset $S \subset \RR \times M, \{W_t(x)\}_{(t,x) \in S}$ is a centered Gaussian random vector;
\itiv $\E[W_s(x) W_t(y)] = (s \wedge t) C_{\mu}(x,y).$
\edes
\blem
\mylabel{th:browncmu} There exists a  Brownian motion on $C(M)$ with covariance $C_{\mu}.$ \elem
\prf Let
\beqarr d_{C_{\mu}}(x,y) &:=&  \sqrt{ C_{\mu}(x,x) - 2 C_{\mu}(x,y) + C_{\mu}(y,y)}\\
&=& \|\nabla\Q_\mu (V_x-V_y)\|_{\Pi(\mu)}\\
&\leq& K \|V_x - V_y\|
\eeqarr
where the last inequality follows from (\ref{sgest2}).
Then $$d_{C_{\mu}}(x,y) \leq   K d(x,y)$$ and the result  follows from Proposition \ref{th:BMexists} and Remark \ref{rem:entrop} in the appendix. \qed

\medskip
We say that a $C(M)$-valued process $Z$ is an {\em Ornstein-Uhlenbeck process of covariance $C_\mu$ and drift $- G_{\mu}$}
if
\beq
\label{eq:OU1}
Z_t = Z_0 - \int_0^t G_{\mu} Z_s ds + W_t
\eeq
where
\bdes
\iti $W$ is a $C(M)$-valued Brownian motion of covariance $C_\mu;$
\itii $Z_0$ is a $C(M)$-valued random variable;
\itiii $W$ and  $Z_0$ are independent.
\edes
Note that we can think of $Z$ as a solution to
 the linear SDE
$$dZ_t=dW_t-G_{\mu}Z_t dt.$$
It follows from section \ref{sec:OU} in the appendix that such a process exists
and defines a Markov process. Furthermore
\bprop
\mylabel{th:OU1}
Under hypothesis \ref{hyp:mu},
\bdes
\iti $(Z_t)$ converges in law toward a $C(M)$-valued random variable  $Z_{\infty};$
\itii $Z_{\infty}$ is {\em Gaussian}, in the sense that for every finite set $S \subset M, \{Z_{\infty}(x)\}_{x \in S}$ is a centered Gaussian random vector;
\itiii Let $\pi^\mu$ denotes the law of $Z_{\infty}.$ Then $\pi^g$ is characterized by its variance
$$\mathsf{Var}(\pi^\mu) : \cM(M) \to \RR,$$ $$m \mapsto  \E( (m Z_{\infty})^2),$$ and for all $m \in \cM,$
\beqarr
\mathsf{Var}(\pi^\mu)(m)  &=&  \int_0^{\infty} \int_{M \times M}
C_{\mu}(x,y) m_t(dx) m_t(dy) dt\\
&=& \int_0^\infty \widehat{C}_\mu(Vm_t,Vm_t) dt \eeqarr
where $$m_t = e^{-t G^*_{\mu}} m.$$
\edes
\eprop
\prf This follows from Proposition \ref{th:weakgauss3} in the appendix. Example \ref{Akernel} shows that assertion $(iii)$ of this proposition is satisfied.
\qed
\subsection{The process $Z^g.$}
For $g=(g_1,\dots,g_n)\in C(M)^n,$ let $\tilde{M} = \{1,\ldots,n\} \cup M $ be the disjoint union of $\{1,\ldots, n\}$ and  $M,$
and $C_{\mu}^g : \tilde{M} \times \tilde{M} \to \RR$ be the function defined by
$$C_{\mu}^g(x,y) =\left \{  \begin{array}{l}
\widehat{C}_{\mu}(g_x,g_y) \mbox{ for } x, y  \in \{1,\ldots, n\},\\  C_{\mu}(x,y) \mbox{ for } x,y \in M,\\
\widehat{C}_{\mu}(V_x,g_y) \mbox{ for } x \in M, y \in \{1, \ldots,
n\}. \end{array} \right.$$
Then $C^g_\mu$ is a Mercer kernel (see section \ref{sec:bmcm}).

A {\em Brownian motion on $\RR^n \times
C(M)$ with covariance $C_{\mu}^g$} is a  $\RR^n  \times C(M)$-valued
stochastic process    $(W^g,W) = \{(W^{g_1}_t,\dots,W^{g_n}_t,
W_t)\}_{t\geq 0}$  such that: \bdes \iti $W = \{W_t\}_{t \geq 0}$ is
a $C(M)$-valued Brownian motion with covariance $C_{\mu};$ \itii For
every finite subset $S \subset \RR \times M, \{ W_t^g,
W_t(x)\}_{(t,x) \in S}$ is a centered Gaussian random vector; \itiii
$\E(W^{g_i}_s W^{g_j}_t) = (s \wedge t) \widehat{C}_{\mu}(g_i,g_j)$
and \\ $\E(W_s(x) W_t^{g_i}) = (s \wedge t)
\widehat{C}_{\mu}(V_x,g_i).$ \edes \blem There  exists a  Brownian
motion on $\RR^n \times C(M)$ with covariance $C_{\mu}^g.$ \elem
\prf Let $\tilde{d}$ be the distance on $\tilde{M}$  defined by
$$\tilde{d}(x,y) =\left \{  \begin{array}{l}
\one_{x \neq y} \mbox{ for } x, y  \in \{1,\ldots, n\},\\
 d(x,y) \mbox{ for } x,y \in M,\\
d(x,x_0) + 1 \mbox{ for } x \in M, y \in \{1, \ldots, n\} \end{array} \right.$$
where $x_0$ is some arbitrary point in $M.$
This makes $\tilde{M}$ a compact metric space, and it is easy to show that the function $C_{\mu}^g$
verifies hypothesis \ref{hyp:entrop} (use the proof of Lemma \ref{th:browncmu}). The result follows by application of Proposition \ref{th:BMexists}. \qed

\medskip
Let now be $Z^g_t=(Z_t^{g_1},\dots,Z_t^{g_n}) \in \RR^n$ denote the  solution to the SDE
\beq
\label{eq:OU1bis}
dZ_t^{g_i} = dW_t^{g_i} - \left(Z^{g_i}_t/2+\hbox{Cov}_{\mu}(Z_t,g_i)\right)dt, \: i = 1, \ldots, n
\eeq
where $(W^g,W)$ is as above and  $Z = (Z_t)$ is given by (\ref{eq:OU1}).

The following result generalizes  Proposition \ref{th:OU1}.
\bprop
\mylabel{th:OU2} Under hypothesis \ref{hyp:mu},
\bdes
\iti  The process $(Z^g_t,Z_t)$  converges in law toward a  centered $\RR^n \times C(M)$ valued Gaussian random variable $(Z^g_{\infty}, Z_{\infty}).$
\itii Let  $\pi^{g,\mu}$ denotes the law of $(Z^g_{\infty}, Z_{\infty}).$ Then $\pi^{g,\mu}$ is characterized by its variance
$$\mathsf{Var}(\pi^{g,\mu}) : \RR^n \times \cM(M) \to \RR,$$
$$(u,m) \mapsto \E \left ( (m Z_{\infty} +  \langle u, Z^g_{\infty}\rangle)^2 \right );$$
and
for all $u \in \RR^n,  m \in \cM(M),$
$$\mathsf{Var} (\pi^{g,\mu})(u,m)=\int_0^\infty \widehat{C}_\mu(f_t,f_t) dt$$
with
$$f_t=e^{-t/2}\sum_i u_ig_i+Vm_t,$$
and where $m_t$ is defined by
\beq\label{eq:mt} m_tf=m_0(e^{-tG_\mu}f) + \sum_{i=1}^n u_i\int_0^t e^{-s/2}\hbox{Cov}_\mu(g_i,e^{-(t-s)G_\mu}f)ds.\eeq
\edes\eprop

\prf Let
$G_\mu^g:\RR^n\times C(M)\to\RR^n\times C(M)$ be the operator
defined  by \beq \label{eq:Gmug}
G_{\mu}^g = \left ( \begin{array}{cc} I/2 &  A_{\mu}^g \\
0 &  G_{\mu} \end{array} \right )
\eeq
 where $A_\mu^g:C(M)\to \RR^n$ is the linear map defined by
$$A_{\mu}^g(f) =
{\Big (}\hbox{Cov}_{\mu}(f,g_1),\dots,\hbox{Cov}_{\mu}(f,g_n)
{\Big )}.$$
Then $(Z^g,Z)$ is a $C(\tilde{M})$-valued Ornstein-Uhlenbeck process of covariance $C_\mu^g$ and drift $-G_\mu^g$.
It is not hard to verify that under hypothesis \ref{hyp:mu}, the assumptions of  Proposition \ref{th:weakgauss3} hold, so that  $(Z^g_t,Z_t)$  converges in law toward a  centered $\RR^n \times C(M)$ valued Gaussian random variable $(Z^g_{\infty}, Z_{\infty})$ with variance
$$\mathsf{Var} (\pi^{g,\mu})(u,m)=\int_0^\infty \widehat{C}_\mu(f_t,f_t) dt$$
with $f_t=\sum_i u_t(i)g_i+Vm_t$ and where $(u_t, m_t) = e^{- t {(G^g_{\mu}})^*}(u,m).$
Now
$$
(G_{\mu}^g)^* = \left ( \begin{array}{cc} I/2 &  0\\
(A_{\mu}^g)^* &  (G_{\mu})^* \end{array} \right )$$
and  $(A_{\mu}^g)^* u = \sum_i u_i (g_i - \mu g_i) \mu.$
 Thus
$u_t = e^{-t/2} u $ and
$$\frac{dm_t}{dt} = - (A_{\mu}^g)^* u_t - (G_{\mu})^* m_t$$
Thus $m_t$ is the solution with $m_0=m$ of
\beq \label{eq:dmt} \frac{dm_t}{dt} = -
e^{-t/2} \left(\sum_i u_i (g_i - \mu g_i)\right) \mu - G^*_{\mu} m_t \eeq
Note that (\ref{eq:dmt}) is equivalent to
$$\frac{d}{dt}(m_t f) = - e^{-t/2} \hbox{Cov}_\mu\left(\sum_i u_i g_i,f\right) - m_t(G_{\mu}f)$$
for all $f\in C(M)$, and $m_0=m$.
From which we deduce that
$$m_t=e^{-tG_\mu^*}m_0-\int_0^t e^{-s/2}e^{-(t-s)G_\mu^*}\left(\sum_i u_i (g_i - \mu g_i)\mu\right) ds$$
which implies the formula for $m_t$ given by (\ref{eq:mt}).
\qed

\medskip
For further reference we call $(Z^g,Z)$ an {\em Ornstein-Uhlenbeck
process of covariance $C_{\mu}^g$ and drift $-G_{\mu}^g$.} It is
called  {\em stationary} when its  initial distribution is
$\pi^{g,\mu}.$
\section{A central limit theorem for $\mu_t$}\label{tlcmu}
We state here the main results of this article.
We assume $\mu^* \in \mathsf{Fix}(\Pi)$ satisfies hypotheses \ref{hyp1} and  \ref{hyp:mu}.
Set $\Delta_t = e^{t/2}(\mu_{e^t}-\mu^*)$, $D_t=V\Delta_t$ and $D_{t+\cdot} = \{D_{t+s}: \, s \geq 0\}.$
Then
\bthm $D_{t+\cdot}$ converges in law, as $t \to \infty,$ towards a stationary Ornstein-Uhlenbeck process of covariance $C_{\mu^*}$ and drift $-G_{\mu^*}.$ \label{mainthm1}\ethm

For $g=(g_1,\dots,g_n)\in C(M)^n$, we set $D^g_t=(\Delta_t g,D_t)$
and $D^g_{t+\cdot} = \{D^g_{t+s}: \, s \geq 0\}.$
Then
\bthm $(D^g_{t+s})_{s\geq 0})$ converges in law towards a
stationary Ornstein-Uhlenbeck process of covariance
$C_{\mu^*}^g$ and drift $-G_{\mu^*}^g$.
\label{mainthm2}\ethm

Define $\widehat{C}:C(M)\times C(M)\to\RR$ the symmetric bilinear form defined by
\beq\label{def:Cchap}
\widehat{C}(f,g)=\int_0^\infty\widehat{C}_{\mu^*}(f_t,g_t)dt, \eeq
with ($g_t$ is defined by the same formula, with $g$ in place of $f$)
\beq\label{def:ftx}f_t(x)=e^{-t/2} f(x)-\int_0^t e^{-s/2}\hbox{Cov}_{\mu^*}(f,e^{-(t-s)G_{\mu^*}}V_x) ds.\eeq
\bcor \label{cor:tcl} $\Delta_t g$ converges in law towards a centered Gaussian variable $Z^g_\infty$ of covariance
$$\E[Z^{g_i}_\infty Z^{g_j}_\infty]=\widehat{C}(f,g).$$
\ecor
\prf Follows from theorem \ref{mainthm2} and the calculus of $\mathsf{Var}(\pi^{g,\mu})(u,0)$. \qed

\subsection{Examples}
\label{secexple}
\subsubsection{Diffusions}
Suppose $V(x,y) = V(x),$ so that $(X_t)$ is just a standard diffusion on $M$
with invariant measure $\mu^* = \frac{exp(-V) \lambda}{\lambda \exp{(-V)}}.$

Let $f\in C(M)$. Then $f_t$ defined by (\ref{def:ftx}) is equal to (using $e^{-tG_{\mu^*}}1=e^{-t/2}1$) $=e^{-t/2} f.$
Thus
\beq\label{cdif} \widehat{C}(f,g)=2\mu^*(f\Q_{\mu^*}g).\eeq
Corollary \ref{cor:tcl} says that
\bthm For all $g\in C(M)^n$,
$\Delta_t^g$ converges in law toward a centered Gaussian variable $(Z_\infty^{g_1},\dots,Z_{\infty}^{g_n})$, with covariance given by
$$\E(Z^{g_i}_{\infty} Z^{g_j}_{\infty}) = 2\mu^*(g_i\Q_{\mu^*}g_j).$$
\ethm
\begin{remark} This central limit theorem for Brownian motions on  compact manifolds has already been considered by Baxter and Brosamler in \cite{bax0} and \cite{bax}; and by Bhattacharya in \cite{bha}  for ergodic diffusions.
\end{remark}

\subsubsection{The case $\mu^*=\lambda$ and $V$ symmetric.}
Suppose here that $\mu^*=\lambda$ and that $V$ is symmetric. We assume (without loss of generality since $\Pi(\lambda)=\lambda$ implies that $V\lambda$ is a constant function) that $V\lambda=0$.


Since $V$ is compact and symmetric, there exists an orthonormal basis $(e_\alpha)_{i\geq 0}$ in $L^2(\lambda)$ and a sequence of reals $(\lambda_\alpha)_{\alpha\geq 0}$ such that $e_0$ is a constant function and
$$V=\sum_{\alpha\geq 1}\lambda_\alpha e_\alpha\otimes e_\alpha.$$
Assume that for all $\alpha$, $1/2+\lambda_\alpha>0$.
Then hypothesis \ref{hyp:mu} holds with $\widehat{\lambda}=\lambda$,
and the convergence of $\mu_t$ towards $\lambda$ holds with positive probability (see \cite{br3}).


Let $f\in C(M)$ and $f_t$ defined by (\ref{def:ftx}), denoting $f^\alpha=\langle f,e_\alpha\rangle_\lambda$ and $f_t^\alpha=\langle f_t,e_\alpha\rangle_\lambda$, we have $f_t^0=e^{-t/2}f^0$ and for $\alpha\ge 1$,
\beqarr
f_t^\alpha
&=& e^{-t/2}f^\alpha - \lambda_\alpha e^{-(1/2+\lambda_\alpha)t}\left(\frac{e^{\lambda_\alpha t}-1}{\lambda_\alpha}\right)f^\alpha\\
&=& e^{-(1/2+\lambda_\alpha)t} f^\alpha.
\eeqarr
Using the fact that
$$\widehat{C}_\lambda(f,g) = 2\lambda (f\Q_\lambda g),$$
this implies that
$$\widehat{C}(f,g) = 2\sum_{\alpha\ge 1}\sum_{\beta\ge 1} \frac{1}{1+\lambda_\alpha+\lambda_\beta}\langle f, e_\alpha\rangle_\lambda \langle g, e_\beta\rangle_\lambda \lambda(e_\alpha\Q_{\lambda}e_\beta).$$

This, with corollary \ref{cor:tcl}, proves
\bthm Assume hypothesis \ref{hyp1} and that $1/2+\lambda_\alpha>0$ for all $\alpha$. Then for all $g\in C(M)^n$,
$\Delta_t^g$ converges in law toward a centered Gaussian variable $(Z_\infty^{g_1},\dots,Z_{\infty}^{g_n})$, with covariance given by
$$\E(Z^{g_i}_{\infty} Z^{g_j}_{\infty}) = \widehat{C}(g_i,g_j).$$
\ethm

In particular,
$$\E(Z^{e_\alpha}_{\infty} Z^{e_\beta}_{\infty}) = \frac{2}{1+\lambda_\alpha+\lambda_\beta} \lambda(e_\alpha\Q_{\lambda}e_\beta).$$
Note that when all $\lambda_\alpha$ are positive, which corresponds to what is named a self-repelling interaction in \cite{br3}, the rate of convergence of $\mu_t$ towards $\lambda$ is bigger than when there is no interaction, and the bigger is the interaction (that is larger $\lambda_\alpha$'s) faster is the convergence.

\section{Proof of the main results} \label{Vmu}
We assume hypothesis \ref{hyp1} and $\mu^*$ satisfies hypothesis \ref{hyp:mu}. It is possible to choose $\kappa$ in hypothesis \ref{hyp:mu} such that $\kappa <1/2$. In the following $\kappa$ will denote such constant. Note that we have $\lambda(-G_{\mu^*})<-\kappa$.
Such $\kappa$ exists when hypothesis \ref{hyp:mu} holds.
\subsection{A lemma satisfied by $\Q_\mu$}
We denote by $\cX(M)$ the space of continuous vector fields on $M$, and equip the spaces $\cP(M)$ and $\cX(M)$ respectively with the weak convergence topology and with the uniform convergence topology.
\blem \label{contgqmu} For all $f\in C(M)$, the mapping $\mu\mapsto \nabla \Q_\mu f$ is a continuous mapping from $\cP(M)$ in $\cX(M)$.
\elem
\prf Let $\mu$ and $\nu$ be in $\cM(M)$, and $f\in C(M)$. Set $g=\Q_\mu f$. Then $f=-A_\mu g+ \Pi(\mu)f$ and
\beqarr
\|\nabla\Q_\mu f -\nabla \Q_\nu f\|_\infty
&=&  \|-\nabla \Q_\mu A_\mu g + \nabla\Q_\nu A_\mu g\|_\infty\\
&=& \|\nabla g +\nabla\Q_\nu A_\mu g\|_\infty\\
&\leq& \|\nabla (g +\Q_\nu A_\nu g)\|_\infty + \|\nabla \Q_\nu (A_\mu - A_\nu)g\|_\infty
\eeqarr
since $\nabla (g +\Q_\nu A_\nu g)=0$ and $(A_\mu - A_\nu)g=\langle \nabla V_{\mu-\nu},\nabla g\rangle$, we get
\beq\label{eq:mu-nu}\|\nabla\Q_\mu f -\nabla \Q_\nu f\|_\infty
\leq K\|\langle \nabla V_{\mu-\nu},\nabla g\rangle\|_\infty.\eeq
Using the fact that $(x,y)\mapsto \nabla V_x(y)$ is uniformly continuous, the right hand term of (\ref{eq:mu-nu}) converges towards 0, when $d(\mu,\nu)$ converges towards 0, $d$ being a distance compatible with the weak convergence. \qed

\subsection{The process $\Delta$}
Set $h_t=V\mu_t$ and $h^*=V\mu^*$.
Recall $\Delta_t=e^{t/2}(\mu_{e^t}-\mu^*)$ and $D_t=V\Delta_t$. Note that $D_t(x)=\Delta_t V_x$.

\medskip
To simplify the notation, we set
$K_s=K_{\mu_s}$, $\Q_s=\Q_{\mu_s}$ and $A_s=A_{\mu_s}$. Let $(M^f_t)_{t\geq 1}$ be the martingale defined by $$M^f_t=\sum_\alpha\int_1^t e_\alpha(\Q_sf)(X_s)dB^\alpha_s.$$
The quadratic covariation of $M^f$ and $M^g$ (with $f$ and $g$ in $C(M)$) is given by
$$\la M^f,M^g\ra_t=\int_1^t \la\nabla \Q_{s}f,\nabla\Q_s g\ra (X_s)ds.$$
Then for all $t\geq 1$ (with $\dot{\Q}_t=\frac{d}{dt}\Q_t$) ,
$$\Q_tf(X_t)-\Q_1f(X_1) = M^f_t + \int_1^t
\dot{\Q}_sf(X_s)ds - \int_1^t K_sf(X_s)ds.$$
Thus
\beqarr \mu_t f &=& \frac{1}{t}\int_1^t K_sf(X_s)ds + \frac{1}{t}\int_1^t \Pi(\mu_s)f ds +\frac{1}{t}\int_0^1 f(X_s)ds\\
&=& -\frac{1}{t}
\left(\Q_tf(X_t)-\Q_1f(X_1)-\int_1^t\dot{\Q}_sf(X_s)ds\right)\\
&&+~\frac{M^f_t}{t} + \frac{1}{t}\int_1^t \la \xi(h_s),f\ra_\lambda ds + \frac{1}{t}\int_0^1 f(X_s) ds.
\eeqarr

Note that $(D_t)$ is a continuous process taking its values in $C(M)$ and that $D_t=e^{t/2}(h_{e^t}-h^*)$.
For $f\in C(M)$ (using the fact that $\mu^*f=\langle \xi(h^*),f\rangle_\lambda$),
\beq \label{eqdelta} \Delta_tf=\sum_{i=1}^5\Delta_t^if \eeq
with
\beqarr
\Delta^1_tf &=& e^{-t/2}\left(-\Q_{e^t}f(X_{e^t})+\Q_1f(X_1) +
\int_1^{e^t} \dot{\Q}_sf(X_s)ds\right)\\
\Delta^2_tf &=& e^{-t/2}M^f_{e^t}\\
\Delta^3_tf &=& e^{-t/2}\int_1^{e^t}\la \xi(h_s)-\xi(h^*)-D\xi(h^*)(h_s-h^*),f\ra_\lambda ds  \\
\Delta^4_tf &=& e^{-t/2}\int_1^{e^t}\la D\xi(h^*)(h_s-h^*),f\ra_\lambda ds  \\
\Delta^5_tf &=& e^{-t/2}\left(\int_0^1 f(X_s)ds -  \mu^* f\right).
\eeqarr

Then $D_t=\sum_{i=1}^5 D^i_t$, where $D^i_t=V\Delta^i_t$.
Finally, note that
\beq \langle D\xi(h^*)(h-h^*),f\rangle_\lambda = -\hbox{Cov}_{\mu^*}(h-h^*,f).\eeq

\subsection{First estimates}
We recall some estimates from \cite{blr}: There exists a constant $K$ such
that for all $f\in C(M)$ and $t> 0$,
\beqarr \|\Q_tf\|_\infty\leq K\|f\|_\infty\\
\|\nabla\Q_tf\|_\infty\leq K\|f\|_\infty\\
\|\dot{\Q}_tf\|_\infty\leq \frac{K}{t}\|f\|_\infty.
\eeqarr
These estimates imply in particular that
$$\langle M^f-M^g\rangle_t\leq K \|f-g\|_\infty \times t$$
and that
\blem \label{lemdom} There exists a constant $K$ depending on
$\|V\|_\infty$ such that for all $t\geq 1$, and all $f\in C(M)$
\beq \label{domineh1} \|\Delta^1_t f\|_\infty+\|\Delta^5_tf\|_\infty\leq K\times (1+t)e^{-t/2}\|f\|_\infty, \eeq
which implies that $((\Delta^1+\Delta^5)_{t+s})_{s\geq 0}$ and $((D^1+D^5)_{t+s})_{s\geq 0}$ both converge towards $0$ (respectively in $\cM(M)$ and in $C(\RR^+\times M)$). \elem

\medskip
We also have
\blem\label{lemd234} There exists a constant $K$ such that for all
$t\geq 0$ and all $f\in C(M)$,
\beqarr
\E[(\Delta^2_t f)^2] &\leq& K\|f\|_\infty^2,\\
|\Delta^3_t f|
 &\leq& K\|f\|_\lambda \times e^{-t/2}\int_0^{t} \|D_s\|^2_\lambda ds,\\
|\Delta^4_t f|
 &\leq& K\|f\|_\lambda\times e^{-t/2}\int_0^{t} e^{s/2}\|D_s\|_\lambda ds.
\eeqarr
\elem
\prf The first estimate follows from
\beqarr
\E[(\Delta^2_tf)^2]
&=& e^{-t}\E[(M^{f}_{e^t})^2] =  e^{-t}\E[\la M^{f}\ra_{e^t}]\\
&\leq& e^{-t}\int_1^{e^t}\|\nabla \Q_s f\|_\infty^2ds \\
&\leq& K \|f\|_\infty^2.
\eeqarr
The second estimate follows from the fact that
$$\| \xi(h)-\xi(h^*)-D\xi(h^*)(h-h^*)\|_\lambda
= O(\|h-h^*\|^2_\lambda).$$
The last estimate
follows easily after having remarked that
\beqarr
|\la D\xi(h^*)(h_s-h^*),f\ra|
&=&|\hbox{Cov}_{\mu^*}(h_s-h^*,f)|\\
&\leq& K\|f\|_\lambda\times\|h_s-h^*\|_{\lambda}\\
&\leq& K\|f\|_\lambda\times s^{-1/2}\|D_{\log(s)}\|_\lambda .
\eeqarr
This proves this lemma. \qed

\subsection{The processes $\Delta'$ and $D'$}
Set $\Delta'=\Delta^2+\Delta^3+\Delta^4$ and $D'=D^2+D^3+D^4$.
For $g\in C(M)$, set
$$\epsilon_t^g = e^{t/2} \la \xi(h_{e^t})-\xi(h^*)-D\xi(h^*)(h_{e^t}-h^*),g\ra_\lambda.$$
Then
$$ d\Delta'_t g = -\frac{\Delta'_t g}{2}dt + dN_t^g +  \epsilon_t^g dt +  \la D\xi(h^*)(D_t),g\ra_\lambda dt $$
where for all $g\in C(M)$, $N^g$ is a martingale. Moreover, for $f$ and $g$ in $C(M)$,
$$\la N^f,N^g\ra_t = \int_0^t \la\nabla \Q_{e^s}f(X_{e^s}),\nabla\Q_{e^s}g(X_{e^s})\ra ds.$$

Then, for all $x$,
$$
dD'_t(x)
= -\frac{D'_t(x)}{2}dt + dM_t(x) +  \epsilon_t(x) dt + \la D\xi(h^*)(D_t),V_x\ra_\lambda dt
$$
where $M$ is the martingale in $C(M)$ defined by $M(x)=N^{V_x}$
and $\epsilon_t(x)=\epsilon_t^{V_x}.$
We also have
$$G_{\mu^*}(D')_t(x) = \frac{D'_t(x)}{2} -\langle D\xi(h^*)(D'_t),V_x\rangle_{\lambda}.$$
Denoting $L_{\mu^*}=L_{-G_{\mu^*}}$ (defined by equation (\ref{defLA}) in the appendix), this implies that
\beqarr
dL_{\mu^*}(D')_t(x)
&=& dD'_t(x)+G_{\mu^*}(D')_t(x)dt\\
&=&  dM_t(x) +\langle D\xi(h^*)((D^1+D^5)_t),V_x\rangle_{\lambda} dt + \epsilon_t(x)dt
\eeqarr
Thus
$$L_{\mu^*}(D')_t(x)=M_t(x)+ \int_0^t\epsilon'_s(x)ds$$
with $\epsilon'_s(x)={\epsilon'}_s V_x$
where for all $f\in C(M)$,
$$\eps'_sf=\epsilon_s^f+\langle D\xi(h^*)((D^1+D^5)_s),f\rangle_\lambda.$$
Using lemma \ref{lemLA},
\beq
D'_t = L_{\mu^*}^{-1}(M)_t + \int_0^t e^{-(t-s)G_{\mu_*}}\epsilon'_sds. \label{eq:D}
\eeq

For $g=(g_1,\dots,g_n)\in C(M)^n$, we denote $\Delta'_t g=(\Delta'_t g_1,\dots,\Delta'_t g_n)$, $N^g=(N^{g_1},\dots,N^{g_n})$ and $\eps'_t g=(\eps'_t g_1,\dots,\eps'_t g_n)$. Then, denoting $L^g_{\mu^*}=L_{-G^g_{\mu^*}}$ (with $G^g_{\mu^*}$ defined by (\ref{eq:Gmug})) we have
$$ L^g_{\mu^*}(\Delta' g , D')_t
=  (N_t^g,M_t) + \int_0^t ({\epsilon'}_sg,\epsilon'_s)ds$$
so that (using lemma \ref{lemLA} and integrating by parts)
\beq(\Delta'_t g,D'_t) = (L_{\mu^*}^g)^{-1}(N^g,M)_t + \int_0^t e^{-(t-s)G^g_{\mu^*}}({\epsilon'}_sg,\epsilon'_s) ds.\label{expD'}\eeq
Moreover
$$ (L_{\mu^*}^g)^{-1}(N^g,M)_t = \left(\widehat{N}^{g_1}_t,\dots,\widehat{N}^{g_n}_t, L_{\mu^*}^{-1}(M)_t\right), $$
where
$$ \widehat{N}^{g_i}_t = N_t^{g_i} - \int_0^t\left(\frac{N^{g_i}_s}{2}
+\widehat{C}_{\mu^*}(L_{\mu^*}^{-1}(M)_s,g_i)\right)ds.$$

\subsection{Estimation of $\epsilon'_t$}
\subsubsection{Estimation of $\|L_{\mu^*}^{-1}(M)_t\|_{\lambda}$}
\blem \label{lem:L-1}
\bdes \iti For all $\alpha\geq 2$, there exists a constant $K_\alpha$ such that for all $t\geq 0$,
$$\E[\|L_{\mu^*}^{-1}(M)_t\|_{\lambda}^\alpha]^{1/\alpha}\leq K_\alpha.$$
\itii a.s. there exists $C$ with $\E[C]<\infty$ such that for all $t\geq 0$,
$$\|L_{\mu^*}^{-1}(M)_t\|_{\lambda} \leq C(1+t).$$
\edes \elem
\prf Since $\|L_{\mu^*}^{-1}(M)_t\|_{\lambda}\le K\|L_{\mu^*}^{-1}(M)_t\|_{\widehat{\lambda}}$, we estimate $\|L_{\mu^*}^{-1}(M)_t\|_{\widehat{\lambda}}$.
We have
$$d L_{\mu^*}^{-1}(M)_t = dM_t - G_{\mu^*}L_{\mu^*}^{-1}(M)_tdt.$$
Let $N$ be the martingale defined by
$$N_t=\int_0^t \left\la
\frac{L_{\mu^*}^{-1}(M)_s}{\|L_{\mu^*}^{-1}(M)_s\|_{\widehat{\lambda}}},dM_s\right\ra_{\widehat{\lambda}}.$$
We have $\la N\ra_t\leq Kt$ for some constant $K$.
Then
\beqarr
d \|L_{\mu^*}^{-1}(M)_t\|^2_{\widehat{\lambda}}
&=& 2 \|L_{\mu^*}^{-1}(M)_t\|_{\widehat{\lambda}} dN_t  -
2\la L_{\mu^*}^{-1}(M)_t,G_{\mu^*}L_{\mu^*}^{-1}(M)_t \ra_{\widehat{\lambda}} dt\\
&& + \quad d\left(\int \la M(x)\ra_t \widehat{\lambda}(dx)\right). \eeqarr
Note that there exists a constant $K$ such that
$$\frac{d}{dt}\left(\int \la M(x)\ra_t
\widehat{\lambda}(dx)\right) \leq K$$
and that (see hypothesis \ref{hyp:mu})
$$\la L_{\mu^*}^{-1}(M)_t,
G_{\mu^*}L_{\mu^*}^{-1}(M)_t \ra_{\widehat{\lambda}}\geq \kappa \|L_{\mu^*}^{-1}(M)_t\|^2_{\widehat{\lambda}}.$$
This implies that
$$\frac{d}{dt}\E[\|L_{\mu^*}^{-1}(M)_t\|^2_{\widehat{\lambda}}]
\leq -2\kappa \E[\|L_{\mu^*}^{-1}(M)_t\|^2_{\widehat{\lambda}}] + K$$
which implies (i) for $\alpha=2$. For $\alpha>2$, we find that
\beqarr
\frac{d}{dt}\E[\|L_{\mu^*}^{-1}(M)_t\|^\alpha_{\widehat{\lambda}}]
&\leq& -\alpha \kappa \E[\|L_{\mu^*}^{-1}(M)_t\|^\alpha_{\widehat{\lambda}}] + K\E[\|L_{\mu^*}^{-1}(M)_t\|^{\alpha-2}_{\widehat{\lambda}}]\\
&\leq& -\alpha \kappa \E[\|L_{\mu^*}^{-1}(M)_t\|^\alpha_{\widehat{\lambda}}] + K\E[\|L_{\mu^*}^{-1}(M)_t\|^\alpha_{\widehat{\lambda}}]^{\frac{\alpha-2}{\alpha}}
\eeqarr
which implies that $\E[\|L_{\mu^*}^{-1}(M)_t\|^\alpha_{\widehat{\lambda}}]$ is bounded.

\medskip
We now prove (ii). Fix $\alpha>1$. Then there exists a constant $K$ such that
$$\frac{\|L_{\mu^*}^{-1}(M)_t\|_{\widehat{\lambda}}^2}{(1+t)^\alpha}
\leq \|L_{\mu^*}^{-1}(M)_0\|_{\widehat{\lambda}}^2 + 2\int_0^t\frac{\|L_{\mu^*}^{-1}(M)_s\|_{\widehat{\lambda}}}{(1+s)^\alpha}dN_s + K. $$
Then BDG inequality implies that
$$\E\left[\sup_{t\geq 0}\frac{\|L_{\mu^*}^{-1}(M)_t\|_{\widehat{\lambda}}^2}{(1+t)^\alpha}\right]
\leq K + 2\sup_{t\geq 0}\left(\int_0^t\frac{Kds}{(1+s)^{2\alpha}}\right)^{1/2}$$
which is finite. This implies the lemma by taking $\alpha=2$. \qed

\subsubsection{Estimation of $\|D_t\|_\lambda$}
Note that $|\eps_t^g|\leq K e^{-t/2}\|D_t\|^2_\lambda\times\|g\|$. Thus
$$|{\eps'}_t g|\leq K e^{-t/2}(1+t+\|D_t\|^2_\lambda)\times \|g\|.$$
This implies (using lemma \ref{lemhyp} and the fact that $0<\kappa<1/2$)
\blem There exists $K$ such that
\beq \label{eps't-s}\left\|\int_0^t e^{-(t-s)G_{\mu^*}} \eps'_s ds\right\|
\leq Ke^{-\kappa t}\left(1+\int_0^t e^{-(1/2-\kappa) s}\|D_s\|_\lambda^2 ds\right).\eeq
\elem


This lemma with lemma \ref{lem:L-1}-(ii) imply the following
\blem a.s. there exists $C$ with $\E[C]<\infty$ such that
\beq \|D_t\|_\lambda\leq  C\times
\left[1+t+\int_0^te^{-s/2}\|D_s\|_\lambda^2ds\right]. \eeq
\elem
\prf  First note that
$$\|D_t\|_\lambda\leq\|D'_t\|_\lambda+K(1+t)e^{-t/2}.$$
Using the expression of $D'_t$ given by (\ref{eq:D}), we get
\beqarr
\|D'_t\|_\lambda
&\leq& \|L_{\mu^*}^{-1}(M)_t\|_\lambda + \left\|\int_0^t e^{-(t-s)G_{\mu^*}}\eps'_s ds\right\|\\
&\leq& C(1+t)+Ke^{-\kappa t}\left(1+\int_0^t e^{-(1/2-\kappa)s} \|D_s\|^2_\lambda ds\right)
\eeqarr
which implies the lemma. \qed

\blem Let $x$ and $\eps$ be real functions. If for all $t\geq 0$,
$$x_t\leq \alpha+\int_0^t\eps_s x_s ds,$$
where $\alpha$ is a real constant, then
$$x_t\leq \alpha\exp\left(\int_0^t\eps_s ds\right).$$
\elem
\prf Similarly to the proof of Gronwall's lemma, we set
$y_t=\int_0^t\eps_s x_s ds$. Then,
$$ \dot{y}_t\leq \alpha\eps_t+\eps_ty_t.$$
Take $\lambda_t=y_t\exp\left(-\int_0^t \eps_s ds\right)$, then
$$\dot{\lambda}_t\leq \alpha\eps_t\exp\left(-\int_0^t \eps_s
ds\right)$$
and
\beqarr
y_t
&\leq& \alpha\int_0^t \eps_s\exp\left(\int_s^t \eps_u
du\right)ds\\
&\leq& \alpha\exp\left(\int_0^t \eps_u du\right) - \alpha.
\eeqarr
This implies the lemma. \qed

\medskip
This lemma implies that $$\|D_t\|_\lambda\leq C(1+t) \times
\exp\left(C\int_0^t
e^{-s/2}\|D_s\|_\lambda ds\right).$$
Since hypothesis \ref{hyp1} implies that
$\lim_{s\to\infty}e^{-s/2}\|D_s\|_\lambda = 0$, this proves that
a.s. for all $\eps>0$, there exists $C_\eps$ such that
$$\|D_t\|_\lambda\leq C_\eps e^{\eps t}.$$
Take $\eps<1/4$. Then
$$\int_0^\infty e^{-s/2}\|D_s\|^2_\lambda ds\leq C_\eps.$$
This implies
\blem a.s., there exists $C$ such that for all $t$,
$$\|D_t\|_\lambda\leq C(1+t).$$ \label{lemestdt} \elem

\subsubsection{Estimation of $\epsilon'_t$}
\blem \label{lem:eps'}a.s. there exists $C$ such that for all $f\in C(M)$,
\beqarr
|{\eps'}_t f|&\leq& C(1+t)^2e^{-t/2} \|f\| 
\eeqarr \elem
\prf We have $|{\eps'}_t f| \leq |\eps_t^f| + K(1+t)e^{-t/2}\|f\|$ and
\beqarr
|\eps_t^f| &\leq & K\|f\|_\lambda\times e^{-t/2}\|D_t\|^2_\lambda \\
&\leq& C\|f\|\times (1+t)^2e^{-t/2} \eeqarr
by lemma \ref{lemestdt}. \qed


\subsection{Estimation of $\|D_t-L_{\mu^*}^{-1}(M)_t\|$}
\blem \label{lem410} $\|D_t-L_{\mu^*}^{-1}(M)_t\|\le C e^{-\kappa t}$.
\elem
\prf We have $\|D_t-D'_t\|\le K(1+t)e^{-t/2}$. So to prove this lemma, it suffices to prove that (see the expression of $D'_t$ given by (\ref{eq:D}))
$$\left\|\int_0^te^{-(t-s)G_{\mu^*}}\eps'_s ds\right\|\le Ce^{-\kappa t}.$$
This term is dominated by
$$K\int_0^t e^{-\kappa(t-s)}\|\eps'_s\| ds.$$
Using the previous lemma, it is also dominated by
$$Ce^{-\kappa t}\int_0^t e^{\kappa s}(1+s)^2e^{-s/2}ds\leq Ce^{-\kappa t}$$
because $\kappa\in ]0,1/2[$. The lemma is proved. \qed

\medskip
In addition, for $g=(g_1,\dots,g_n)\in C(M)^n$, setting $$\Delta_t g=(\Delta_t g_1,\dots,\Delta_t g_n),$$
\blem $\|(\Delta_t g,D_t)-(L^g_{\mu^*})^{-1}(N^g,M)_t\|\le C(1+\|g\|)e^{-\kappa t}$.
\label{lem411}\elem
\prf We have $\|(\Delta_t g,D_t)-(\Delta'_t g,D'_t)\|\leq K(1+\|g\|)(1+t)e^{-\kappa t}$. So to prove this lemma, using (\ref{expD'}), it suffices to prove that
\beq \label{taest}\left\|\int_0^te^{-(t-s)G^g_{\mu^*}}({\eps'}_s g,\eps'_s) ds\right\| \leq K(1+\|g\|)e^{-\kappa t}. \eeq
Using hypothesis \ref{hyp:mu} and the definition of $G^g_{\mu^*}$, we have that for all positive $t$,
$$\|e^{-t G^g_{\mu^*}}\|\leq Ke^{-\kappa t}.$$
This implies
$$\|e^{-(t-s)G^g_{\mu^*}}({\eps'}_s g,\eps'_s)\|\leq Ke^{-\kappa(t-s)}\|\eps'_s\|\times (1+\|g\|).$$
Thus the term (\ref{taest}) is dominated by
$$K(1+\|g\|)\int_0^t e^{-\kappa(t-s)}\|\eps'_s\| ds,$$
from which we prove (\ref{taest}) like in the previous lemma. \qed

\subsection{Tightness results}

We refer the reader to section \ref{sec:tighness} in the appendix, where tightness criteria for families of $C(M)$-valued random variables are given. They will be used in this section.

\subsubsection{Tightness of $(L_{\mu^*}^{-1}(M)_t)_{t\geq 0}$}
In this section we prove the following lemma which in particular implies the tightness of $(D_t)_{t\geq 0}$ and of $(D'_t)_{t\geq 0}$.
\blem \label{esttens}$(L_{\mu^*}^{-1}(M)_t)_{t\geq 0}$ is tight. \elem

\prf
We have the relation (that defines $L_{\mu^*}^{-1}(M)$)
$$dL_{\mu^*}^{-1}(M)_t(x)=-G_{\mu^*}L_{\mu^*}^{-1}(M)_t(x) dt + dM_t(x).$$
Thus, using the expression of $G_{\mu^*}$
$$dL_{\mu^*}^{-1}(M)_t(x)=-\frac{1}{2}L_{\mu^*}^{-1}(M)_t(x) dt + A_t(x)dt + dM_t(x),$$
with
$$A_t(x)=\widehat{C}_{\mu^*}(V_x,L_{\mu^*}^{-1}(M)_t).$$
Since $\mu^*$ is absolutely continuous with respect to $\lambda$, we have that
$$\|A_t\|\leq K\|L_{\mu^*}^{-1}(M)_t\|_\lambda$$
and therefore (using lemma \ref{lem:L-1} (i) for $\alpha=2$)
$$\sup_t\E[\|A_t\|^2]<\infty.$$
We also have
$$\hbox{Lip}(A_t)\leq K\|L^{-1}_{\mu^*}(M)_t\|_\lambda,$$
where $\hbox{Lip}(A_t)$ is the Lipschitz constant of $A_t$ (see (\ref{deflip}))

In order to prove this tightness result, we first prove that for all $x$, $(L_{\mu^*}^{-1}(M)_t(x))_t$ is tight. Setting $Z_t^x=L_{\mu^*}^{-1}(M)_t(x)$ we have
\beqarr
\frac{d}{dt}\E[(Z_t^x)^2]
&\leq& -\E[(Z_t^x)^2] + 2\E[|Z_t^x|\times |A_t(x)|] + \frac{d}{dt}\E[\langle M(x)\rangle_t] \\
&\leq&  -\E[(Z_t^x)^2] + K \E[(Z_t^x))^2]^{1/2} + K
\eeqarr
which implies that $(L_{\mu^*}^{-1}(M)_t(x))_t$ is bounded in $L^2(\P)$ and thus tight.

We now estimate $\E[|Z_t^x-Z_t^y|^\alpha]^{1/\alpha}$ for $\alpha$ greater than $2$ and the dimension of $M$. Setting $Z_t^{x,y}=Z_t^x-Z_t^y$, we have
\beqarr
\frac{d}{dt}\E[(Z_t^{x,y})^\alpha]
&\leq& -\frac{\alpha}{2}\E[(Z_t^{x,y})^\alpha] + \alpha\E[(Z_t^{x,y})^{\alpha-1}|A_t(x)-A_t(y)|] \\
&& + \frac{\alpha(\alpha-1)}{2}\E\left[(Z_t^{x,y})^{\alpha-2}\frac{d}{dt}\langle M(x)-M(y)\rangle_t\right] \\
&\leq& -\frac{\alpha}{2}\E[(Z_t^{x,y})^\alpha] + \alpha d(x,y)\E[(Z_t^{x,y})^{\alpha-1}\hbox{Lip}(A_t)] \\
&&+ Kd(x,y)^2\E[(Z_t^{x,y})^{\alpha-2}] \\
&\leq& -\frac{\alpha}{2}\E[(Z_t^{x,y})^\alpha] + Kd(x,y)\E[(Z_t^{x,y})^{\alpha-1}\|L^{-1}(M)_t\|_\lambda] \\
&& + Kd(x,y)^2 \E[(Z_t^{x,y})^{\alpha-2}] \\
&\leq& -\frac{\alpha}{2}\E[(Z_t^{x,y})^\alpha] + Kd(x,y)\E[(Z_t^{x,y})^\alpha]^{\frac{\alpha-1}{\alpha}}\E[\|L^{-1}(M)_t\|^\alpha_\lambda]^{1/\alpha} \\
&&+ Kd(x,y)^2\E[(Z_t^{x,y})^\alpha]^{\frac{\alpha-2}{\alpha}} \\
&\leq& -\frac{\alpha}{2}\E[(Z_t^{x,y})^\alpha] + Kd(x,y)\E[(Z_t^{x,y})^\alpha]^{\frac{\alpha-1}{\alpha}} \\
&&+ Kd(x,y)^2\E[(Z_t^{x,y})^\alpha]^{\frac{\alpha-2}{\alpha}}.
\eeqarr
Thus, if $x_t=\E[(Z_t^{x,y})^\alpha]/d(x,y)^\alpha$,
$$\frac{dx_t}dt\le -\frac{\alpha}{2}x_t + Kx_t^{\frac{\alpha-1}{\alpha}} + Kx_t^{\frac{\alpha-2}{\alpha}}.$$
It is now an exercise to show that $x_t\le K$ and so that $$\E[(Z_t^{x,y})^\alpha]^{1/\alpha}\leq Kd(x,y).$$ Using corollary \ref{tightcrit}, this completes the proof for the tightness of $(L_{\mu^*}^{-1}(M)_t)_t$. \qed

\brem \label{rqK} Kolmogorov's theorem (see theorem 1.4.1 and its proof in Kunita (1990)), with the estimates given in the proof of this lemma, implies that $$\sup_t\E[\|L_{\mu^*}^{-1}(M)_t\|]<\infty.$$
\erem

\subsubsection{Tightness of $((L_{\mu^*}^g)^{-1}(N^g,M)_t)_{t\geq 0}$}
Fix $g=(g_1,\dots,g_n)\in C(M)^n$.
Let $\widehat{\Delta}g$ be defined by the relation
$$(\widehat{\Delta}g,L_{\mu^*}^{-1}(M))=(L_{\mu^*}^g)^{-1}(N^g,M).$$
Set $A_tg=(A_tg_1,\dots,A_tg_n)$ with
$A_tg_i=\widehat{C}_{\mu^*}(g_i,L_{\mu^*}^{-1}(M)_t)$. Then
$$d\widehat{\Delta}_tg=dN_t^g-\frac{\widehat{\Delta}_tg}{2} dt + A_tg dt.$$
Thus,
$$\widehat{\Delta}_tg=e^{-t/2}\int_0^t e^{s/2}dN_s^g + e^{-t/2}\int_0^t e^{s/2}A_sg ds.$$
Using this expression it is easy to prove that $(\widehat{\Delta}_tg)_{t\geq 0}$ is bounded in $L^2(\P)$. This implies, using also lemma \ref{esttens}
\blem \label{esttensg}$((L^g_{\mu^*})^{-1}(N^g,M)_t)_{t\geq 0}$ is tight. \elem

\subsection{Convergence in law of $(N^g,M)_{t+\cdot}-(N^g,M)_t$}
In this section, we denote by $\E_t$ the conditional expectation with respect to $\cF_{e^t}$. We also set $\Q=\Q_{\mu^*}$ and $C=\widehat{C}_{\mu^*}$.

\subsubsection{Preliminary lemmas.}
For $f\in C(M)$ and $t\geq 0$, set $N^{f,t}_s=N^f_{t+s}-N^f_t$.
\blem \label{convcov} For all $f$ and $g$ in $C(M)$,
$$\lim_{t\to\infty} \langle N^{f,t},N^{g,t}\rangle_s=s\times C(f,g).$$
\elem
\prf
Set $$G(z)=\langle\nabla \Q f,\nabla \Q g\rangle(z) -C(f,g)$$ and
$$G_u(z)=\langle\nabla \Q_u f,\nabla \Q_u g\rangle(z) -C(f,g).$$
We have
\beqarr
\la N^{f,t},N^{g,t}\ra_s -s\times C(f,g) &=&
\int_{e^t}^{e^{t+s}} G_u(X_u) \frac{du}{u}\\
&=& \int_{e^t}^{e^{t+s}} (G_u-G)(X_u) \frac{du}{u}\\
&&+\quad\int_{e^t}^{e^{t+s}}G(X_u) \frac{du}{u}
.
\eeqarr

Integrating by parts, we get that
$$\int_{e^t}^{e^{t+s}}G(X_u) \frac{du}{u} = (\mu_{e^{t+s}}G-\mu_{e^t}G) + \int_0^s (\mu_{e^{t+u}} G)du. $$
Since $\mu^* G=0$, this converges towards $0$ on the event $\{\mu_t\to\mu^*\}$. The
term $\int_{e^t}^{e^{t+s}} (G_u-G)(X_u) \frac{du}{u}$ converges towards $0$ because $(\mu,z)\mapsto \nabla
\Q_\mu f(z)$ is continuous. This proves the lemma. \qed

\medskip
Let $f_1,\dots,f_n$ be in $C(M)$.
Let $(t_k)$ be an increasing sequence converging to $\infty$ such that
the conditional law of $M^{n,k}=(N^{f_1,t_k},\dots,N^{f_n,t_k})$ given $\cF_{e^{t_k}}$ converges in law towards a $\RR^n$-valued process
$W^n=(W_1,\dots,W_n)$.

\blem \label{convlfd} $W^n$ is a centered Gaussian process such that for all $i$ and $j$,
$$\E[W^n_i(s)W^n_j(t)]=(s\wedge t)  C(f_i,f_j).$$ \elem
\prf We first prove that $W^n$ is a martingale.
For all $k$, $M^{n,k}$ is a martingale.
For all $u\leq v$, B\"urkholder-Davies-Gundy inequality (BDG inequality in the following) implies that
$(M^{n,k}(v)-M^{n,k}(u))_k$ is bounded in $L^2$.

Let $l\geq 1$, $\p\in C(\RR^l)$, $0\leq s_1\leq \cdots\leq s_l\leq u$
and $(i_1,\dots,i_l)\in\{1,\dots,n\}^l$.
Then for all $k$ and $i\in\{1,\dots,n\}$, the martingale property implies that
$$\E_{t_k}[(M^{n,k}_i(v)-M^{n,k}_i(u))Z_k]=0$$
where $Z_k$ is of the form
\beq\label{Zk}Z_k=\p(M^{n,k}_{i_1}(s_1),\dots,M^{n,k}_{i_l}({s_l})).\eeq
Using the convergence of the conditional law of $M^{n,k}$ given $\cF_{e^{t_k}}$ towards the law of
$W^n$ and since $(M_i^{n,k}(v)-M^{n,k}_i(u))_k$ is
uniformly integrable (because it is bounded in $L^2$), we prove that
$$\E[(W^n_i(v)-W^n_i(u))Z]=0$$
where  $Z$ is of the form
\beq\label{Z}Z=\p(W^n_{i_1}(s_1),\dots,W^n_{i_l}(s_l)).\eeq
This implies that $W^n$ is a martingale.

\medskip
We now prove that for $(i,j)\in\{1,\dots,n\}$ (with $C=C_{\mu^*}$),
$$\la W^n_i,W^n_j\rangle_s=s\times C(f_i,f_j).$$
By definition of $\la M^{n,k}_i,M^{n,k}_j\rangle$ (in the following
$\langle\cdot,\cdot\rangle^v_u =
\langle\cdot,\cdot\rangle_v-\langle\cdot,\cdot\rangle_u$)
\beqar\label{mnk}
\E_{t_k}\left[\left((M^{n,k}_i(v)-M^{n,k}_i(u))(M^{n,k}_j(v)-M^{n,k}_j(u))\right.\right.&&\\
&&\hskip-80pt \left.\left.-\la M^{n,k}_i,M^{n,k}_j\ra^v_u\right) Z_k\right]=0\non \eeqar
where $Z_k$ is of the form (\ref{Zk}).
Using the convergence in law and the fact that
$(M^{n,k}(v)-M^{n,k}(u))^2_k$ is bounded in $L^2$ (still using BDG inequality),
we prove that as $k\to\infty$,
$$ \E_{t_k}[(M^{n,k}_i(v)-M^{n,k}_i(u))(M^{n,k}_j(v)-M^{n,k}_j(u))Z_k] $$
converges towards
$$ \E[(W^{n}_i(v)-W^{n}_i(u))(W^{n}_j(v)-W^{n}_j(u))Z] .$$
with $Z$ of the form (\ref{Z}). Now,
\beqarr
\E_{t_k}[\la M^{n,k}_i,M^{n,k}_j\ra_v Z_k]
-v \times\E[Z]\times C(x_i,x_j)
&&\\
&&\hskip-120pt=\quad \E_{t_k}[(\la M^{n,k}_i,M^{n,k}_j\ra_v-v\times C(f_i,f_j))Z_k]\\
&&\hskip-100pt+\quad v\times (\E_{t_k}[Z_k]-\E[Z])\times C(f_i,f_j)
\eeqarr
The convergence in $L^2$ of $\la M^{n,k}_i,M^{n,k}_j\ra_v$ towards
$v\times C(f_i,f_j)$ shows that the first term converges towards
$0$. The convergence of the conditional law of $M^{n,k}$ with respect to $\cF_{e^{t_k}}$ towards $W^n$  shows that
the second term converges towards $0$. Thus
$$\E\left[\left((W^n_i(v)-W^n_i(u))(W^n_j(v)-W^n_j(u))
  -(v-u) C(f_i,f_j)\right) Z\right]=0.$$
This shows that $\la W^n_i,W^n_j\ra_s=s\times C(f_i,f_j)$.
We conclude using L\'evy's theorem. \qed

\subsubsection{Convergence in law of $M_{t+\cdot}-M_t$}
In this section, we denote by $\cL_t$ the conditional law of
$M_{t+\cdot}-M_t$ knowing $\cF_{e^t}$. Then $\cL_t$ is a probability
measure on $C(\RR^+\times M)$.

\bprop When $t\to\infty$, $\cL_t$ converges weakly towards the law of
a $C(M)$-valued Brownian motion of covariance $C_{\mu^*}$. \eprop
\prf In the following, we will simply denote $M_{t+\cdot}-M_t$ by
$M^t$. We first prove that
\blem $\{\cL_t:~t\geq 0\}$ is tight. \elem
\prf For all $x\in M$, $t$ and $u$ in $\RR^+$,
\beqarr
\E_t[(M^t_u(x))^2]
&=& \E_t\left[\int_t^{t+u} d\langle M(x)\rangle_s\right]\\
&\leq& Ku.
\eeqarr
This implies that for all $u\in\RR^+$ and $x\in M$, $(M^t_u(x))_{t\geq 0}$ is tight.

\smallskip
Let $\alpha>0$. We fix $T>0$.
Then for $(u,x)$ and $(v,y)$ in $[0,T]\times M$,
using BDG inequality,
\beqarr
\E_t[|M^t_u(x)-M^t_v(y)|^\alpha]^{\frac{1}{\alpha}} &\leq& \E_t[|M^t_u(x)-M^t_u(y)|^\alpha]^{\frac{1}{\alpha}}\\
&&+\quad \E_t[|M^t_u(y)-M^t_v(y)|^\alpha]^{\frac{1}{\alpha}} \\
&\leq& K_\alpha \times(\sqrt{T}d(x,y)+\sqrt{|v-u|})
\eeqarr
where $K_\alpha$ is a positive constant depending only on $\alpha$,
$\|V\|$ and $\hbox{Lip}(V)$ the Lipschitz constant of $V$.

We now let $D_T$ be the distance on
$[0,T]\times M$ defined by
$$D_T((u,x),(v,y))=K_\alpha \times
(\sqrt{T}d(x,y)+\sqrt{|v-u|}).$$
The covering number $N([0,T]\times M,D_T,\eps)$ is of order $\eps^{-d-1/2}$ as $\eps\to 0$.
Taking $\alpha>d+1/2$, we conclude using corollary \ref{tightcrit}. \qed

\medskip
Let $(t_k)$ be an increasing sequence converging to $\infty$
and $N$ a $C(M)$-valued random process (or a $C(\RR^+\times M)$ random
variable) such that $\cL_{t_k}$ converges in law towards $N$.

\blem $N$ is a $C(M)$-valued Brownian motion of covariance $C_{\mu^*}$. \elem
\prf Let $W$ be a $C(M)$-valued Brownian motion of covariance $C_{\mu^*}$.
Using lemma \ref{convlfd}, we prove that for all $(x_1,\dots,x_n)\in M^n$, $(N(x_1),\dots,N(x_n))$ has the same distribution as $(W(x_1),\dots,X(x_n))$. This implies the lemma. \qed

\medskip
Since $\{\cL_t\}$ is tight, this lemma implies that $\cL_t$ converges weakly towards the law of a $C(M)$-valued Brownian motion of covariance $C_{\mu^*}$. \qed

\subsubsection{Convergence in law of $(N^g,M)_{t+\cdot}-(N^g,M)_t$}
In this section, we fix $g=(g_1,\dots,g_n)\in C(M)^n$ and we denote by $\cL^g_t$ the conditional law of
$(N^g,M)_{t+\cdot}-(N^g,M)_t$ knowing $\cF_{e^t}$. Then $\cL^g_t$ is a probability
measure on $C(\RR^+\times M\cup\{1,\dots,n\})$. In the following we will denote $(N^{g,t},M^t)$ the process $(N^g,M)_{t+\cdot}-(N^g,M)_t$.

Let $(W_t^f)_{(t,f)\in \RR^+\times C(M)}$ be a $\mathcal{X}(M)$-valued Brownian motion of covariance $\widehat{C}_{\mu^*}$. Denoting $W_t(x)=W_t^{V_x}$, then $W=(W_t(x))_{(t,x)\in\RR^+\times M}$ is a $C(M)$-valued Brownian motion of covariance $C_{\mu^*}$. For $g=(g_1,\dots,g_n)\in C(M)^n$, $W^g$ will denote $(W^{g_1},\dots,W^{g_n})$. In the following we will simply denote $(W^g,W)$ the process $(W_t^g,(W_t(x))_{x\in M})_{t\geq 0}$.

\bprop As $t$ goes to $\infty$, $\cL^g_t$ converges weakly towards the law of $(W^g,W)$. \eprop
\prf  We first prove that
\blem $\{\cL^g_t:~t\geq 0\}$ is tight. \elem
\prf This is a straightforward consequence of the tightness of $\{\cL_t\}$ and of the fact that
for all $\alpha>0$, there exists $K_\alpha$ such that for all nonnegative $u$ and $v$, $\E_t[|N^{g,t}_u-N^{g,t}_v|^\alpha]^{\frac{1}{\alpha}} \leq K_\alpha \sqrt{|v-u|}.$
\qed

\medskip Let $(t_k)$ be an increasing sequence converging to $\infty$
and $(\tilde{N}^g,\tilde{M})$ a $\RR^n\times C(M)$-valued random process (or a $C(\RR^+\times M\cup\{1,\dots,n\})$ random
variable) such that $\cL^g_{t_k}$ converges in law towards $(\tilde{N}^g,\tilde{M})$.
Then lemmas \ref{convcov} and \ref{convlfd} imply that $(\tilde{N}^g,\tilde{M})$ has the same law as $(W^g,W)$.
Since $\{\cL^g_t\}$ is tight, $\cL^g_t$ convergences towards the law of $(W^g,W)$. \qed

\subsection{Convergence in law of $D$}
\subsubsection{Convergence in law of $(D_{t+s}-e^{-sG_{\mu^*}}D_t)_{s\geq 0}$}\label{sec491}
We have
\beqarr
D'_{t+s}-e^{-sG_{\mu^*}}D'_t &=&  L_{\mu^*}^{-1}(M^t)_s +
\int_0^s e^{-(s-u)G_{\mu^*}}\eps'_{t+u} du.
\eeqarr
Since (using lemma \ref{lem:eps'})
$$ \left\|\int_0^s e^{-(s-u)G_{\mu^*}} \eps'_{t+u} du\right\|
\leq K e^{-\kappa t} $$
and $\|D_t-D'_t\|\leq K(1+t)e^{-t/2}$,
this proves that
$(D_{t+s}-e^{-sG_{\mu^*}} D_t - L_{\mu^*}^{-1}(M_{t+\cdot}-M_t)_s)_{s\geq 0}$ converges towards
$0$. Since $L_{\mu^*}^{-1}$ is continuous, this proves that the law of
$L_{\mu^*}^{-1}(M_{t+\cdot}-M_t)$ converges weakly
towards $L_{\mu^*}^{-1}(W)$. Since $L_{\mu^*}^{-1}(W)$ is an Ornstein-Uhlenbeck process of covariance $C_{\mu^*}$ and drift $-G_{\mu^*}$ started from $0$, we have
\bthm The conditional law of $(D_{t+s}-e^{-sG_{\mu^*}}D_t)_{s\geq 0}$ given $\cF_{e^t}$ converges weakly
towards an Ornstein-Uhlenbeck process of covariance $C_{\mu^*}$ and drift $-G_{\mu^*}$ started from $0$. \ethm

\subsubsection{Convergence in law of $D_{t+\cdot}$}\label{sec492}
We can now prove theorem \ref{mainthm1}. We here denote by $\P_t$ the semigroup of an Ornstein-Uhlenbeck process of covariance $C_{\mu^*}$ and drift $-G_{\mu^*}$, and we denote by $\pi$ its invariant probability measure.

\medskip
We know that (as $t\to\infty$) $(D_{t+s}-e^{-s G_{\mu^*}}D_t)_{s\geq 0}$ converges in law towards $L_{\mu^*}^{-1}(W)$, where $W$ is a $C(M)$-valued Brownian motion of covariance $C_{\mu^*}$.
Since $(D_t)_{t\geq 0}$ is tight, there exists $\nu\in\cP(C(M))$ and an increasing sequence $t_n$ converging towards $\infty$ such that $D_{t_n}$ converges in law towards $\nu$. Then $D_{t_n+\cdot}$ converges in law towards $(L_{\mu^*}^{-1}(W)_s+e^{-s G_{\mu^*}}Z_0)$, with $Z_0$ independent of $W$ and distributed like $\nu$. This proves that $D_{t_n+\cdot}$ converges in law towards an Ornstein-Uhlenbeck process of covariance $C_{\mu^*}$ and drift $-G_{\mu^*}$.

We now fix $t>0$. Let $s_n$ be a
subsequence of $t_n$ such that $D_{s_n-t+\cdot}$ converges in law. Then
$D_{s_n-t}$ converges towards a law we denote by $\nu_t$ and $D_{s_n-t+\cdot}$ converges in law towards an Ornstein-Uhlenbeck process of covariance $C_{\mu^*}$ and drift $-G_{\mu^*}$.
Since $D_{s_n}=D_{s_n-t+t}$, $D_{s_n}$ converges in law towards $\nu_t\P_t$.
On the other hand $D_{s_n}$ converges in law towards $\nu$. Thus $\nu_t\P_t=\nu$.

Let $\p$ be a Lipschitz bounded function on $C(M)$. Then
\beqarr
|\nu_t\P_t\p-\pi \p|
&=& \left|\int (\P_t\p(f)-\pi\p)\nu_t(df)\right|\\
&\leq& \int |\P_t\p(f)-\P_t\p(0)|\nu_t(df) + |\P_t\p(0)-\pi\p|
\eeqarr
where the second term converges towards $0$ (using (\ref{eq:convptp})) and the first term is dominated by (using lemma \ref{lem:ptpf-g})
$$Ke^{-\kappa t} \int \|f\| \nu_t(df).$$

It is easy to check that
\beqarr
\int \|f\| \nu_t(df)
&=& \lim_{k\to\infty}\int (\|f\|\wedge k) \nu_t(df)\\
&=& \lim_{k\to\infty} \lim_{n\to\infty}\E[\|D_{s_n-t}\|\wedge k]\\
&\leq& \lim_{n\to\infty}\E[\|D_{s_n-t}\|]\\
&\leq& \sup_t\E[\|D_t\|].
\eeqarr
Since
\beqarr
\|D_t\|
&\leq& \|D^1_t+D^5_t\| + \|L^{-1}_{\mu^*}(M)_t\|\\
&& + \left\|\int_0^t e^{(t-s)G_{\mu^*}}\eps'_s ds\right\|,
\eeqarr
using the estimates (\ref{domineh1}), the proof of lemma \ref{lem410} and remark \ref{rqK}, we get that
$$\sup_{t\geq 0}\E[\|D_t\|]<\infty.$$

Taking the limit, we prove $\nu\p=\pi\p$ for all Lipschitz bounded function $\p$ on $C(M)$.
This implies $\nu=\pi$, which proves the theorem. \qed

\subsubsection{Convergence in law of $D^g$}
We can also prove theorem \ref{mainthm2}.

\medskip
For $g=(g_1,\dots,g_n)\in C(M)^n$, we set $D^g_t=(\Delta_t g,D_t)$, and ${D'_t}^g=(\Delta'_t g,D'_t)$.
Since $\|D_t^g-{D'}_t^g\|\leq K(1+t)e^{-t/2}$, instead of studying $D^g$, we can only study ${D'_t}^g$.
Then
\begin{eqnarray*}
{D'}_{t+s}^g-e^{-sG^g_{\mu^*}}{D'}_t^g
&=& (L^g_{\mu^*})^{-1}(N^{g,t},M^t)_s\\
&& + \int_0^s e^{-(s-u)G^g_{\mu^*}}(\eps'_{t+u}g,\eps'_{t+u})du.
\end{eqnarray*}
The norm of the second term of the right hand side (using the proof of lemma \ref{lem411}) is dominated by
\begin{eqnarray*}
&\leq& K(1+\|g\|)\int_0^s e^{-\kappa(s-u)} \|\eps'_{t+u}\| du\\
&\leq& K\int_0^s e^{-\kappa(s-u)}(1+t+u)^2e^{-(t+u)/2}du \\
&\leq& ke^{-\kappa t}
\end{eqnarray*}
Like in section \ref{sec491}, since $(L^g_{\mu^*})^{-1}(W^g,W)$ is an Ornstein-Uhlenbeck process of covariance $C^g_{\mu^*}$ and drift $-G^g_{\mu^*}$ started from $0$,
\bthm The conditional law of $((\Delta^g,D)_{t+s}-e^{-sG^g_{\mu^*}}(\Delta^g,D)_t)_{s\geq 0}$ given $\cF_{e^t}$ converges weakly towards an Ornstein-Uhlenbeck process of covariance $C^g_{\mu^*}$ and drift $-G^g_{\mu^*}$ started from $0$. \ethm

From this theorem, like in section \ref{sec492}, we prove theorem \ref{mainthm2}. \qed

\section{Appendix : Random variables and Ornstein-Uhlenbeck processes on $C(M)$}
\mylabel{sec:app}
\subsection{$C(M)$-valued random variables}
\subsubsection{Generalities}

Let $(M,d)$ be a compact metric space (note that there is no
assumption here that $M$ is a manifold), $C(M)$ the space of real
valued continuous functions on $M$ equipped with the uniform norm $\|f\| =
\sup_{x \in M} |f(x)|.$ By classical results, $C(M)$ is a separable
(see e.g~ \cite{stroock}) Banach space (see e.g~ \cite{dudley} or
\cite{conway}) and its topological dual is the space $\cM(M)$ of
bounded signed measures on $M$ (see e.g~\cite{dudley} or
\cite{conway}).  For $\mu \in \cM(M)$ and $f \in C(M)$ we use the
notation $\mu f = \langle \mu, f \rangle = \int_M f d\mu$.

Let $(\Omega, \cF, \P)$ be a probability space. A {\em $C(M)$-valued random variable} is a Borel map $F :\Omega \to C(M).$

For $x \in M,$ let  $\pi_x : C(M) \to \RR,$ denote the projection defined by
$$\pi_x(f) = f(x).$$
\blem
\mylabel{th:generalities}
 The Borel $\sigma$-field on $C(M)$ is the $\sigma$-field generated by the maps $\{\pi_x\}_{ x \in M}$. In particular
\bdes
\iti A map $F : \Omega \to C(M)$ is a $C(M)$-valued random variable if and only if $\{\pi_x(F)\}_{x \in M}$ is a family of  real valued random variables.
\itii The law of a $C(M)$-valued random variable is determined by its finite dimensional distributions (i.e~ the law of $\{\pi_x(F)\}_{x \in I}$ with $I \subset M$ finite).
\edes
\elem
\prf
Let $\cA = \sigma\{\pi_x, x \in M\}$ and $\cB$ the Borel $\sigma$-field on $C(M).$
The maps $\pi_x$ being continuous, $\cB$ contains $\cA.$
Conversely, let $B_f(r) = \{g \in C(M) \: : \|g-f\|\leq r\}$ and let $S$ be a countable dense subset of $M.$ Then $B_f(r) = \cap_{x \in S} \{g \in C(M) : \: |\pi_x(f) - \pi_x(g)| \leq r\}$ Hence $B_f(r) \in \cA.$ Since $C(M)$ is separable, $\cB$ is generated by  the sets $\{B_f(r), f \in C(M), r \geq 0 \}.$
\qed
\subsubsection{Tightness criteria}\label{sec:tighness}
 Let $\cP(C(M))$ be the space of Borel probability measures on
 $C(M)$.
An element $\nu$ of $\cP(C(M))$ is the law of a
$C(M)$-valued random variable $F$, and $\nu=\P_F$. Recall that a sequence
$\{\nu^n\}$ in $\cP(C(M))$ is said {\em converging weakly} towards $\nu\in\cP(C(M))$ if $\int \p d\nu_n \to \int \p d\nu$ for every bounded and continuous function $\p : C(M) \to  \RR.$
A sequence $\{F_n\}$ of $C(M)$-valued random variable is said {\em converging in law} towards $F$ a $C(M)$-valued random variable if $\{P_{F_n}\}$ converges in law towards $\P_F$.
A family $\cX \subset \cP(C(M))$  is said to be {\em tight} if for every $\eps > 0$ there exists some compact set $\cK \subset C(M)$ such that $\P(\cK) \geq 1-\eps$ for all $\cP \in \cX.$ A family of random variables is said to be tight if the family of their laws is tight.

Since $C(M)$ is a  separable and complete,  Prohorov theorem \cite{billingsley} asserts that  $\cX \subset \cP(C(M))$ is tight if and only if it is relatively compact.

The next proposition gives a useful criterium for a class of random
variables to be tight. It follows directly from  \cite{lt}
(Corollary 11.7 p. 307 and the remark following Theorem 11.2). A
function $\psi : \Rp \to \Rp$ is a Young function if it is convex,
increasing and $\psi(0) = 0.$ If $Z$ is a real valued random
variable, we let
$$\|Z\|_{\psi} = \inf \{c > 0 \: : \E\big(\psi(|Z|/c)\big) \leq 1 \}.$$
 For $\eps>0$, we denote by $N(M,d;\eps)$ the covering number of $E$ by balls of radius less than $\eps$ (i.e. the minimal number of balls of radius less than $\eps$ that cover $E$), and by $D$ the diameter of $M$.
\bprop \mylabel{th:tightcrit}
Let $(F_t)_{t\in I}$ be a family
of $C(M)$-valued random variables and $\psi$ a Young function.
Assume that
\bdes\iti There exists $x\in E$ such that $(F_t(x))_{t\in I}$ is tight;
\itii $\|F_t(x) - F_t(y)\|_{\psi} \leq K d(x,y);$
\itiii $\int_0^D \psi^{-1}(N(M,d;\eps))d \eps <\infty.$
\edes
Then $(F_t)_{t\geq 0}$ is tight. \eprop

\medskip
\bcor
\mylabel{tightcrit}
 Suppose $M$ is a compact finite dimensional  manifold of dimension $r,$  $d$ the Riemannian distance,   and
$$[\E|F_t(x)-F_t(y)|^{\alpha}]^{1/\alpha} \leq K d(x,y)$$ for some $\alpha > r.$ Then conditions $(ii)$ and $(iii)$ of Proposition \ref{th:tightcrit} hold true.
\ecor
\prf  One has
$N(E,d;\eps)$ is of order $\eps^{-r};$ and for  $\psi(x) = x^{\alpha}, \|\cdot\|_{\psi}$ is the $L^{\alpha}$ norm. Hence the result.
\qed

\subsubsection{$C(M)$-valued Gaussian variable}

Recall that a (centered)  real-valued random variable $Y$ with variance $\sigma^2$ is said to be {\em  Gaussian} if it has distribution
$$P_{Y}(dx) = \frac{1}{\sqrt{2\pi} \sigma} \exp(-\frac{x^2}{2\sigma^2}) dx.$$
Its characteristic function is then
$$\Phi_Y(t) = \E[\exp(it Y)] = \exp(-\frac{t^2 \sigma^2}{2}).$$
Here we adopt the convention that the zero function ($Y = 0$) is Gaussian with variance $0$ and that all the Gaussian random variables are centered.

A family  $\{Y_i\}_{i \in I}$ of real-valued random variables is said to be Gaussian if for all finite set
 $J \subset I$ and for all $\alpha \in \RR^J,$
$\sum_{j \in J} \alpha_j Y_j$ is Gaussian.

A $C(M)$-valued random variable $F$ is said to be {\em Gaussian} if
for all $\mu \in \cM(M)$, $\langle \mu , F \rangle$ is Gaussian.

\blem A $C(M)$-valued random variable $F$ is Gaussian if and only if the family
 $\{\pi_x(F)\}$ is Gaussian.
\elem
\prf
The direct implication is obvious. We prove the second.
Assume that  $\{\pi_x(F)\}$ is a  Gaussian family.  Let $\mu$ be a
 probability over $M.$ By the strong law of large number and the
 separability of $C(M)$ there exists a nonempty set $\Lambda \subset M^{\NN}$
 (actually $\Lambda$ has  $\mu^{\NN}$ measure $1$) such that  for all
 $(x_i) \in \Lambda$ and all $f \in C(M)$
$$\lim_{n \to \infty} \frac{1}{n} \sum_{i = 1}^n f(x_i) = \langle
\mu, f \rangle.$$
In particular
 $Y_n \to   \langle
\mu, F \rangle$
where $$Y_n = \frac{1}{n} \sum_{i = 1}^n F(x_i).$$
And, by Lebesgue theorem, $\Phi_{Y_n}(t) \to \Phi_{\langle \mu, F \rangle}(t).$
Since, by assumption $Y_n$ is Gaussian, $\Phi_{Y_n}(t) = \exp(-t^2
\sigma_n^2/2).$ Let $\sigma \in [0, \infty]$ be a limit point of
$(\sigma_n).$ Then $\Phi_{\langle \mu, F \rangle}(t) = \exp(-t^2 \sigma^2/2).$  This proves that $\sigma < \infty$  (a characteristic
function being continuous) and that $\langle \mu, F \rangle$ is
Gaussian.

If now $\mu \in \cM(M)$ by Jordan-Hann decomposition we may write
$\mu = a \mu_1 - b \mu_2$ with $a,b \geq 0$ and  $\mu_1, \mu_2$
probabilities. It follows from what precede that $\langle \mu,
X\rangle$ is Gaussian. \qed

\medskip
Given a $C(M)$-valued Gaussian random variable $F$ we let
 $\mathsf{Var}_{F} : \cM(M) \to \RR$ denote the {\em variance function} of $F$ defined by $$\mathsf{Var}_F(\mu) = \E(\langle \mu, F \rangle^2).$$
In view of lemma \ref{th:generalities} $(ii),$ the law of $F$ is entirely determined by its variance function.

A useful property of Gaussian variables is the following.
\blem
\label{lem:invarGausslin}
Let $M'$ be another compact metric space and $A : C(M) \to C(M')$ a bounded linear operator. Let $F$ be a $C(M)$-valued Gaussian random variable. Then
$AF$ is a $C(M')$-valued Gaussian random variable with variance
$$\mathsf{Var}_{AF} = \mathsf{Var}_F \circ A^*$$
where $A^* : \cM(M') \to \cM(M)$ is the adjoint of $A.$ \elem \prf
follows from the duality $\langle \mu , AF \rangle = \langle A^*
\mu, F \rangle$ and the definitions. \qed

\subsection{Brownian motions on $C(M)$.} \label{sec:bmcm}
Let $C : M \times M \to \RR$  be  a  continuous symmetric
(i.e~$C(x,y) = C(y,x)$) function such that
$$\sum_{ij} a_i a_j C(x_i,x_j) \geq 0 $$
for every finite sequence $(a_i,x_i)$ with $a_i \in \RR$ and $x_i \in M.$
Such a function is sometimes called a {\em Mercer kernel}.

A  {\em Brownian motion on $C(M)$ with covariance $C$} is a $C(M)$-valued stochastic process  $W  = \{ W_t\}_{ t \geq 0}$ such that $W_0 = 0$ and  for each $T  \geq 0,\, W^T = \{W_t(x): \: t \leq T, x \in M\}$ is a $C([0,T] \times M)$-valued Gaussian random variable with variance
$$\mathsf{Var}_{W^T}(\nu) = \int_{([0,T] \times M)^2} (s \wedge t) C(x,y) \nu(ds dx) \nu(dt dy)$$ or equivalently
$$\E(W_t(x) W_s (y)) = (s \wedge t) C(x,y).$$

Let $$d_C(x,y) = \sqrt{ C(x,x) - 2 C(x,y) + C(y,y)}.$$
The $d_C$ is a pseudo-distance on $M$.
 For $\eps > 0,$ let
 $$\omega_C(\eps) = \sup \{\eta > 0 : d(x,y) \leq \eta \Rightarrow d_C(x,y) \leq \eps\}.$$
Then $N(M,d;\omega_C(\eps)) \geq N(M, d_C; \eps)$.
\begin{hypothesis}
\mylabel{hyp:entrop}
$$\int_0^1 \log (N(d , M; \omega_C(\eps))~ d \eps < \infty$$
where $N(d ,M;\eta)$  is the covering number of $M$ by balls of radius less than $\eta.$
\end{hypothesis}
\begin{remark}
\mylabel{rem:entrop}
Assume that $M$ is a compact finite dimensional manifold and that
$d_C(x,y) \leq K d(x,y)^{\alpha}$ for some $\alpha > 0.$ Then $\omega_C(\eps) \leq (\frac{\eps}{K})^{1/\alpha}$ and $N(d,M;\eta) = O(\eta^{-dim(M)});$ so that the preceding hypothesis holds.
\end{remark}
\bprop
\mylabel{th:BMexists}
Under hypothesis \ref{hyp:entrop}  there exists a Brownian motion on $C(M)$ with covariance $C.$
\eprop
\prf
By Mercer Theorem (see e.g~ \cite{dieudonne}) there exists a countable family of function $\Psi_i \in C(M), i \in \NN,$
such that
$$C(x,y) = \sum_i \Psi_i(x) \Psi_i(y)$$ and the convergence is uniform.
Let $B^i, i \in \NN,$ be a family of independent standard Brownian motions. Set   $$W^n_t(x) = \sum_ {i \leq n}  B^i_t \Psi_i(x), \, n \geq 0.$$
Then, for each $(t,x) \in \RR^{+} \times M,$ the sequence
$(W^n_t(x))_{n \geq 1}$ is a Martingale. It is furthermore
bounded in $L^2$ since $$\E[(W^n_t(x))^2] = t \sum_{i \leq n } \Psi_i(x)^2  \leq t C(x,x).$$
Hence by Doob's convergence theorem one may define
$$W_t(x) = \sum_{i \geq 0}  B^i_t \Psi_i(x).$$
Let now  $S \subset \RR^+ \times M$ be a countable and dense set. It is easily checked that the
 family $(W_t(x))_{(t,x) \in S}$ is a centered Gaussian family with covariance given by
$$\E[W_s(x)W_t(y)]=(s\wedge t) C(x,y),$$
In particular, for $ t \geq s$
\beqarr
\E[(W_s(x) -W_t(y))^2] &=& s C(x,x) - 2 s  C(x,y) + t C(y,y)\\
&=& s d_C(x,y)^2 +  (t-s) C(y,y) \\
&\leq& K (t-s) + s d_C(x,y)^2
\eeqarr
This later bound combined with  classical results on Gaussian processes (see e.g~Theorem 11.17 in \cite{lt}) implies that
$(t,x) \mapsto W_t(x)$ admits a  version uniformly continuous over $S_T = \{(t,x) \in S: \: t \leq T\}$. By density it can be extended to a continuous
 (in $(t,x)$) process
$$W=(W_t(x))_{\{(t,x)\in\RR^+\times M\}}$$
The process $W$ can be viewed as a $C(M)$-valued continuous random
process with the desired covariance.
\qed

\subsection{Ornstein-Ulhenbeck processes}
\mylabel{sec:OU}
Let $A : C(M) \to C(M)$ be a bounded operator and $W$ a
$C(M)$-valued Brownian motion with covariance $C$ as defined in the
preceding section.

An {\em Ornstein-Ulhenbeck process} with drift $A,$ covariance $C$ and initial condition $F_0 = f \in C(M)$ is defined to be a
$C(M)$ valued stochastic process continuous in $t,$ such that
\beq
\label{GenOU}
F_t - f= \int_0^t A F_s ds + W_t.
\eeq
Note that we may think of $F$ as the
solution to the ``stochastic differential equation'' on $C(M):$  $$
dF_t = A F_t dt + dW_t$$ with initial condition  $F_0 = f \in C(M).$

Our aim here is to construct such a solution and state some of its properties.

We let $(e^{tA})_{t \in \RR}$ denote the linear flow induced by $A.$
Recall that for each $t$, $e^{tA}$ is the bounded operator on $C(M)$ defined by
$$e^{tA} = \sum_{k \in \NN} \frac{(tA)^k}{k!}.$$
Given $T > 0$ we let $L_A : C(\RR^+ \times M) \to C(\RR^+ \times M)$
be defined by
\beq\label{defLA} L_A(f)_t = f_t - f_0 - \int_0^t A f_s ds, \qquad t\geq 0.\eeq
Given $T > 0$ we let $L^T_A : C([0,T] \times M) \to C([0,T] \times M)$
be defined by
\beq\label{defLAT} L^T_A(f)_t = f_t - f_0 - \int_0^t A f_s ds, \qquad 0\leq t\leq T.\eeq
Note that if for $f\in C(\RR^+ \times M)$, we let $f^T\in C([0,T] \times M)$ be defined by $f^T_t=f_t$ for $t\in [0,T]$.
\blem\label{lemLAT}
$L_A^T$ is a bounded operator and its restriction to $C_0([0,T] \times
M) = \{f \in C([0,T] \times M) \: : f_0 = 0\}$ is bijective with
inverse $(L^T_A)^{-1}$
defined by
\beq
\label{defLT-1}
(L_A^T)^{-1}(g)_t = g_t + \int_0^t e^{(t-s)A} A g_s ds, \qquad t\in [0,T].
\eeq
\elem
\prf
Linearity of $L_A^T$ is obvious. Also
$$\|L_A^T(f)\| \leq (2 + T\|A\|)  \|f\|.$$
This proves that $L_A^T$ is bounded.

Observe that $L_A^T(f) = 0$ implies that $f_t = e^{tA} f_0.$ Hence $L_A^T$
restricted to $C_0([0,T] \times M)$ is injective.
Let $g \in C_0([0,T] \times M)$ and let $f_t$ be given by the left
hand side of (\ref{defLT-1}). Then
$$h_t = L_A^T(f)_t - g_t = \int_0^t e^{(t-s)A}A g_s ds - \int_0^t A f_s ds.$$
It is easily seen that $h$ is  differentiable and that   $\frac{d}{dt} h_t = 0.$ This proves
that $h_t = h_0 = 0.$
\qed

\medskip
We also have
\blem\label{lemLA}
The restriction of $L_A$ to $C_0(\RR^+ \times M) = \{f \in C(\RR^+ \times M) \: : f_0 = 0\}$ is bijective with
inverse $(L_A)^{-1}$ defined by
\beq
\label{defL-1}
L_A^{-1}(g)_t = g_t + \int_0^t e^{(t-s)A} A g_s ds.
\eeq
\elem

The next lemma easily follows.
\blem
\mylabel{th:gensolOU}
For all $f \in C(M)$ and $g \in C_0(\RR^+ \times M)$  the solution
to $$f_t = f + \int_0^t A f_s ds + g_t,$$
is given by
$$f_t = e^{tA} f + L_A^{-1}(g)_t.$$
\elem

If now $W$ is a $C(M)$-valued Brownian motion as defined in the preceding section, one may define
$$F_t = e^{tA} f + L_A^{-1}(W)_t.$$
Such a process is the unique solution to (\ref{GenOU}). Note that,
by Lemma \ref{lem:invarGausslin} $(F_t - e^{tA} f)_{t \leq T}$ is a
$C_0([0,T] \times M)$-valued Gaussian random variable. In particular
\bprop \mylabel{th:varOU} Let $(F_t)$ be the solution to
(\ref{GenOU}) with initial condition $F_0 = 0.$ Then for each $t
\geq 0$, $F_t$ is a $C(M)$-valued Gaussian random variable with
variance
$$\mathsf{Var}_{F_t}(\mu) = \int_0^t \langle   \mu, e^{sA} C e^{sA^*} \mu \rangle ds.$$
where $C : \cM(M) \to C(M)$ is the operator
defined by $C\mu (x) = \int_M C(x,y) \mu(dy).$
\eprop
\prf
Fix $T > 0.$ To shorten notation let $G : C_0([0,T] \times M) \to C(M)$ be the operator defined by
$$G(g) = (L_A^T)^{-1}(g)_T.$$
Hence $F_T = G(W^T).$
By Lemma  \ref{lem:invarGausslin}, $F_T$ is Gaussian with variance
$$\mathsf{Var}_{F_T} = \mathsf{Var}_{W^T} \circ G^*.$$
Now, for all $\nu \in \cM([0,T] \times M)$
$$\mathsf{Var}_{W^T} (\nu) = \langle \nu,  \cC \nu \rangle$$
where
$$\cC \nu (s,x) = \int_{[0,T]\times M} (s \wedge u) C(x,y) \nu(du dy).$$
Thus
$$\mathsf{Var}_{F_T}(\mu) = \langle \mu,  G \cC G^* \mu \rangle.$$
Our next goal is to compute  $G \cC G^* \mu.$
It easily follows from the definition of $G$ that
$$G^* \mu =  \delta_T \otimes \mu + ds A^* e^{(T-s)A^*} \mu.$$
Thus (integrating by parts)
$$\cC G^* \mu = \cC \nu_1 + \cC \nu_2$$
with $\nu_1 = \delta_T \otimes \mu$
and $\nu_2 = ds A^* e^{(T-s)A*} \mu.$
One has $$\cC \nu_1(s,x) = s (C\mu) (x);$$
$$\cC \nu_2 (s,x) = \int_0^T (s \wedge u) (C \dot{m_u}(\mu))(x) du$$
with
$m_u(\mu) = - e^{(T-u)A^*} \mu$ and $\dot{m_u}(\mu)$ stands for the derivative of $u \mapsto m_u(\mu).$
Thus
$$\cC \nu_2(s,x) = - s C \mu - \int_0^s C m_u(\mu) (x) du$$
and
$$\cC G^* \mu (s,x) = - \int_0^s C m_u(\mu) (x) du.$$
Set  $$h_s(x)  = \int_0^s C m_u(\mu)(x) du.$$
Then
$$G \cC G^* \mu  = h_T + \int_0^T e^{(T-s)A} A h_s ds = \int_0^T e^{(T-s)A} \dot{h}_s ds$$
$$ = \int_0^T e^{(T-s) A} C e^{(T-s)A^*} \mu ds = \int_0^T e^{s A} C e^{s A^*} \mu ds.$$
\qed

\subsubsection{Asymptotic Behaviour}
Let $\lambda(A) = \lim_{t \to \infty} \frac{\log(\|e^{tA}\|}{t}$
which exists by subadditivity. Then for some constant $K<\infty$, $\|e^{tA}\|\le Ke^{\lambda(A)t}$ for all positive $t$.
Let $(F_t)$ denote the solution to (\ref{GenOU}), with $F_0=f\in C(M)$.
\bcor
\mylabel{th:weakgauss1}
Assume $\lambda(A) < 0.$ Then for each $\mu \in \cM(M)$
$(\langle \mu, F_t \rangle)$ converges in law toward a Gaussian random variable with variance
$$\mathsf{V}(\mu) = \int_0^{\infty}  \langle \mu,  e^{sA} C e^{sA*} \mu \rangle ds .$$
\ecor
\prf follows from proposition \ref{th:varOU} and Lemma  \ref{th:gensolOU} \qed
\bcor
\mylabel{th:weakgauss2}  Assume that  $\lambda(A) < 0.$  Set
 $$d_{\mathsf{V}}(x,y) = \sqrt{\mathsf{V}(\delta_x - \delta_y)}$$
and
$$\omega_{\mathsf{V}}(\eps) = \sup \{\eta > 0 : d(x,y) \leq \eta \Rightarrow d_{\mathsf{V}}(x,y) \leq \eps\}.$$   Assume furthermore that
$\omega_{\mathsf{V}}$ verifies the condition expressed by hypothesis
\ref{hyp:entrop}. Then $(F_t)$ converges in law toward a
$C(M)$-valued Gaussian random variable with variance $\mathsf{V}.$
\ecor \prf Let $\nu_t$ denote the law of $F_t.$ Corollary
\ref{th:weakgauss1} and lemma \ref{th:generalities} imply that every
limit point of $\{\nu_t\}$ (for the weak* topology) is the law of a
$C(M)$-valued Gaussian variable  with variance $\mathsf{V}.$ The
proof then reduces to show that $(\nu_t)$ is relatively compact or
equivalently that $\{F_t\}$ is tight. We use  Proposition
\ref{th:tightcrit}. The first condition follows  from Lemma
\ref{th:weakgauss1}. Let $\psi(x) = e^{x^2} -1.$ It is easily
verified that for any real valued Gaussian random variable $Z$ with
variance $\sigma^2$, $\|Z\|_{\Psi} = \sigma \sqrt{8/3}.$ Hence
$\|F_t(x)-F_t(y)\|_{\psi} \leq 2 d_{\mathsf{V}}(x,y)$ so that
condition $(ii)$ holds with the pseudo distance $d_{\mathsf{V}}.$ By
definition of $\omega_{\mathsf{V}},$
$N(M,d;\omega_{\mathsf{V}}(\eps)) \geq N(M, d_{\mathsf{V}}; \eps)$
and since $\psi^{-1}(u) = \sqrt{\log(u-1)}$ condition $(iii)$ is
verified. \qed

\medskip
Denote by $\P_t$ the semigroup associated to an Ornstein-Uhlenbeck process of covariance $C$ and drift $A$. Then for all bounded measurable $\p:C(M)\to\RR$ and $f\in C(M)$,
\beq \label{def:ptou}
\P_t\p(f)=\E[\p(F_t)].\eeq
Denote by $\pi$ the law of a $C(M)$-valued Gaussian random variable with variance $\mathsf{V}.$ Then $\pi$ is the invariant probability measure of $\P_t$, i.e. $\pi\P_t=\pi$. Corollary \ref{th:weakgauss2} implies that, when $\lambda(A)<0$, for all $f\in C(M)$ and all bounded continuous $\p:C(M)\to\RR$,
\beq \label{eq:convptp} \lim_{t\to\infty}\P_t\p(f)=\pi \p.\eeq
Even thought we don't have the speed of convergence in the previous limit, we have
\blem\label{lem:ptpf-g}
Assume that $\lambda(A)<0$. For all bounded Lipschitz continuous $\p:C(M)\to\RR$, all $f$ and $g$ in $C(M)$,
$$|\P_t\p(f)-\P_t\p(g)|\leq Ke^{\lambda(A)t}\|f-g\|.$$
\elem
\prf We have $\P_t\p(f)=\E[\p(L_A^{-1}(W)_t+e^{tA}f)]$. So, using the fact that $\p$ is Lipschitz,
\beqarr
|\P_t\p(f)-\P_t\p(g)| &\le& K\|e^{tA}(f-g)\|\\
&\le& Ke^{\lambda(A)t}\|f-g\|.
\eeqarr
This proves the lemma. \qed
\medskip

To conclude this section we give  a set of simple sufficient conditions  ensuring that the hypotheses of corollary \ref{th:weakgauss2} are satisfied.

For $f \in C(M)$ we let \beq \hbox{Lip}(f) = \sup_{x \neq y} \frac{|f(x)-f(y)|}{d(x,y)} \in \RR^+ \cup \{\infty\}.\label{deflip}\eeq
A map $f$ is said to be Lipschitz provided $Lip(f) < \infty.$
\bprop
\mylabel{th:weakgauss3}
Assume
\bdes
\iti $N(d,M;\eps) = O(\eps^{-r})$ for some $r > 0$ (This holds in particular if $M$ is a finite dimensional manifold).
\itii $x \mapsto C(x,y)$ is Lipschitz uniformly in $y.$ That is
$$\sup_{z\in M}|C(x,z) - C(y,z)| \leq K d(x,y)$$ for some $K \geq 0.$
\itiii There exists $K > 0$ such that
$$Lip(Af) \leq K(Lip(f) + \|f\|).$$
\itiv $\lambda(A) < 0$
\edes
Then the hypotheses, hence the conclusion,  of corollary \ref{th:weakgauss2} are satisfied.
\eprop

We begin with the following lemma.
\blem Under hypotheses $(iii)$ and $(iv)$ of proposition \ref{th:weakgauss3}
$$Lip(e^{tA} f) \leq e^{Kt} (Lip(f) + K' \|f\|)$$
for some constants $K,K'.$
\elem
\prf
For all $x,y$
$$|e^{tA} f(x) - e^{tA} f (y)| = \left|\int_0^t [A e^{sA} f (x) - A e^{sA} f (y)] ds + f(x) - f(y)\right|$$
$$
\leq K \left ( \int_0^t \left [ Lip(e^{sA} f) + \|e^{sA} f\| \right ] ds + Lip(f) \right ) d(x,y).$$
Since $\lambda = \lambda(A) < 0,$ there exists $C>0$ such that
$\|e^{sA}\| \leq Ce^{-s\lambda}.$ Thus
$$Lip(e^{tA}f ) \leq K \int_0^t Lip(e^{sA} f ) ds + \frac{K C}{\lambda}\|f\| + Lip(f)$$
and the result follows from Gronwal's lemma.  \qed

\medskip
We now pass to the proof of the proposition. In what follows the constants  may change from line to line.

\prf Set $\mu = \delta_x - \delta_y$ and  $f_s = C e^{sA^*} \mu$ so
that
$$ \langle \mu, e^{sA} C e^{sA^*} \mu \rangle = e^{sA} f_s(x) - e^{sA} f_s(y).$$
It follows from  hypotheses $(ii)$ and $(iv)$ that $$Lip(f_s) + \|f_s\|  \leq K e^{-s\lambda}$$
for some positive constants $K$ and $a.$
Therefore, by the preceding lemma,
$$Lip(e^{sA} f_s) \leq K e^{s \alpha}$$ for some (other) positive constants $K, \alpha.$
Thus \beqarr d_{\mathsf{V}}(x,y)^2 &\leq& d(x,y) \int_0^T Lip(e^{sA}
f_s) ds + \int_{T}^{\infty} (e^{sA} f(x) - e^{sA} f(y)) ds\\
&\leq& d(x,y) \int_0^T K e^{s \alpha} ds + 2 \int_{T}^{\infty}
\|e^{sA} f_s\|ds\\
&\leq& K \left(d(x,y)e^{\alpha T} + \int_{T}^{\infty} e^{-s\lambda}
ds\right)\\
&\leq& K (d(x,y) e^{\alpha T} + e^{- \lambda T}). \eeqarr Let
$\gamma = \frac{\alpha}{\lambda}, \, \eps > 0, $ and   $T =
-\ln(\eps)/\lambda.$   Then $$d_{\mathsf{V}}^2(x,y) \leq K
(\eps^{-\gamma} d(x,y) + \eps).$$ Therefore
$$d(x,y) \leq \eps^{\gamma+1} \Rightarrow d^2_{\mathsf{V}}(x,y) \leq K \eps,$$ so that
$N(d,M;\omega_ {\mathsf{V}}( \eps))=O\big(\eps^{-2r(\gamma+1)}\big)$
and hypothesis \ref{hyp:entrop} holds true. \qed

\bex
\mylabel{Akernel}
{\rm
Let $$Af(x) = \int f(y) k(x,dy)$$
with
$$k(x,dy) = k_0(x,y) \mu(dy) + \sum_{i = 1}^n a_i(x) \delta_{b_i(x)}$$
where
\bdes
\iti $\mu$ is a bounded measure on $M,$
\itii $k_0(x,y)$ is bounded and uniformly Lipschitz in $x,$
\itiii $a_i : M  \to  \RR$ and  $b_i : M \to M$ are  Lipschitz.
\edes
Then hypothesis $(iii)$ of proposition \ref{th:weakgauss3} is verified. }
\eex

\end{document}